%% file: torus_21_08_05_for_public.tex
\documentclass[reqno]{amsart}
\usepackage{amsmath}
\usepackage{amsthm}
\usepackage{amssymb}
\usepackage{amscd}
\usepackage{amsfonts}
\usepackage{mathrsfs,amsxtra,upref}
\usepackage[mathscr]{eucal}

\numberwithin{equation}{section}

\makeatletter

\def\th@plain{%
\let\thm@indent\noindent
\thm@headfont{\caps}
\let\thmhead\thmhead@plain
\let\swappedhead\swappedhead@plain
\thm@preskip.5\baselineskip\@plus.2\baselineskip\@minus.2\baselineskip
\thm@postskip\thm@preskip
\slshape
}

\def\th@remark{%
\let\thm@indent\noindent
\thm@headfont{\caps}
\let\thmhead\thmhead@plain
\let\swappedhead\swappedhead@plain
\thm@preskip.5\baselineskip\@plus.2\baselineskip\@minus.2\baselineskip
\thm@postskip\thm@preskip
\upshape
}

\def\th@boldremark{%
\let\thm@indent\noindent
\thm@headfont{\bfseries}
\let\thmhead\thmhead@plain
\let\swappedhead\swappedhead@plain
\thm@preskip.5\baselineskip\@plus.2\baselineskip\@minus.2\baselineskip
\thm@postskip\thm@preskip
\upshape
}

\makeatother

\theoremstyle{plain}
\newtheorem{Theorem}{Theorem}[section]
\newtheorem{Theorem--Definition}[Theorem]{Theorem--Definition}
\newtheorem{Corollary}[Theorem]{Corollary}
\newtheorem{Lemma}[Theorem]{Lemma}
\newtheorem{Lemma--Definition}[Theorem]{Lemma--Definition}
\newtheorem{Proposition}[Theorem]{Proposition}
\newtheorem{Proposition--Definition}[Theorem]{Proposition--Definition}
\theoremstyle{remark}
\newtheorem{Remark}[Theorem]{Remark}%
\newtheorem{Remark--Definition}[Theorem]{Remark--Definition}%
\newtheorem{Definition}[Theorem]{Definition}

\newtheorem{Notation}[Theorem]{Notation}

\DeclareMathOperator{\at}{\text{@}}
\DeclareMathOperator{\Arg}{Arg}
\DeclareMathOperator{\Aut}{{Aut}\,}
\DeclareMathOperator{\card}{\#}
\DeclareMathOperator*{\const}{const}

\DeclareMathOperator{\Id}{Id}
\DeclareMathOperator{\id}{id}
\DeclareMathOperator{\Img}{Im}
\DeclareMathOperator{\Ker}{Ker}
\DeclareMathOperator{\supp}{supp}
\DeclareMathOperator{\thh}{{\!}^{th}}

\font\caps=cmcsc10

\font\twelvegtc=eufm10 scaled 1200

\font\ninegtc=eufm9
\font\sevengtc=eufm7

\newfam\gtcfam

\textfont\gtcfam=\twelvegtc
\scriptfont\gtcfam=\ninegtc
\scriptscriptfont\gtcfam=\sevengtc

\def\modo#1{\left| #1 \right|}
\def\Def{\overset{\text{def}}{=\!=}}
\input cd
\LaTeXfeatures

\begin{document}
\title{Configuration spaces of tori}
\author[Y. Feler]{Yoel Feler}
\address{Technion\\ Haifa\\   Israel 32000}
\email{feler$\at$techunix.technion.ac.il}
\thanks{This material is based upon PhD thesis supported by Technion}

\begin{abstract}
\noindent
The configuration spaces
${\mathcal C}^n({\mathbb T}^2)
=\{Q\subset {\mathbb T}^2\,|\ \#Q=n\}$
and 
${\mathcal E}^n({\mathbb T}^2)=\{(q_1,...,q_n)
\in ({\mathbb T}^2)^n\,|\ q_i\ne q_j \ \forall\,i\ne j\}$
of a torus ${\mathbb T}^2
={\mathbb C}/\{\text{\rm a lattice}\}$
are complex manifolds. We prove
that for $n>4$ any holomorphic self-map $F$ of
${\mathcal C}^n({\mathbb T}^2)$ either carries
the whole of ${\mathcal C}^n({\mathbb T}^2)$ into
an orbit of the diagonal $\Aut({\mathbb T}^2)$ action
in ${\mathcal C}^n({\mathbb T}^2)$ or is of the form
$F(Q)=T(Q)Q$, where 
$T\colon{\mathcal C}^n({\mathbb T}^2)\to\Aut({\mathbb T}^2)$
is a holomorphic map. We also prove that for $n>4$
any endomorphism of the torus braid group
$B_n({\mathbb T}^2)=\pi_1({\mathcal C}^n({\mathbb T}^2))$
with a non-abelian image preserves 
the pure torus braid group $PB_n({\mathbb T}^2)
=\pi_1({\mathcal E}^n({\mathbb T}^2))$.
\end{abstract}

\maketitle

\tableofcontents

\section{Introduction}

\noindent In this paper we study certain properties of
the torus configuration spaces and torus braid groups.
\vskip7pt

\subsection{Configuration spaces}\label{ss 1.1}
The $n\thh$ configuration space
${\mathcal C}^n={\mathcal C}^n(X)$ of a complex space $X$ consists
of all $n$ point subsets (``configurations")
$Q=\{q_1,...,q_n\}\subset X$. An automorphism $T\in\Aut X$
gives rise to the holomorphic endomorphism $F_T$ of ${\mathcal C}^n$
defined by $F_T(Q)=TQ=\{Tq_1,...,Tq_n\}$.
If $\Aut X$ is a complex Lie group,
one may take $T=T(Q)$ depending analytically on a configuration
$Q\subset X$ and define the corresponding
holomorphic endomorphism $F_T$ by $F_T(Q)=T(Q)Q$.
Such a map $F_T$ is called {\em tame}. On the other hand,
choosing a base configuration $Q^0=\{q_1^0,...,q_n^0\}\subset X$,
one may define an endomorphism $F_{T,Q^0}$, \
$F_{T,Q^0}(Q)=T(Q)Q^0=\{T(Q)q_1^0,...,T(Q)q_n^0\}$, which
certainly maps the whole configuration space
into one orbit $(\Aut X)Q^0$ of the diagonal $\Aut X$ action in
${\mathcal C}^n$; endomorphisms that have the latter
property are said to be {\em orbit-like}.
\vskip5pt

\noindent Of course, configuration spaces of
a certain specific space $X$ may admit ``sporadic"
holomorphic endomorphisms, which are neither tame nor orbit-like;
but it is not to be expected
that there is a
general construction of such endomorphisms.
\vskip5pt

\noindent Let us restrict ourselves to the simplest
interesting case when $X$ is a non-hyperbolic Riemann
surface\footnote{The automorphism group
of a hyperbolic Riemann surface $X$ is discrete;
in fact, a ``generic" hyperbolic Riemann surface does not admit
a non-identical automorphism.},
i.e. one of the following curves:
complex projective line ${\mathbb{CP}}^1$;
complex affine line $\mathbb C$; complex affine line with one puncture
${\mathbb C}^*=\mathbb C\setminus\{0\}$;
a torus ${\mathbb T}^2=\mathbb{C}/$\{a lattice of rank 2\}.
\vskip5pt

\noindent V. Lin (see \cite{Lin72b,Lin79,LinSphere},
etc.) proved that for $n>4$ and $X=\mathbb C$ or
$X={\mathbb{CP}}^1$ every holomorphic endomorphism $F$
of ${\mathcal C}^n(X)$ is either tame or orbit-like. Moreover,
$F$ is tame if and only if it is non-abelian, meaning that
the image $F_*(\pi_1({\mathcal C}^n(X)))$
of the induced endomorphism $F_*$ of the fundamental group
$\pi_1({\mathcal C}^n(X))$ {\rm (which is the braid group $B_n(X)$
of $X$)} is a non-abelian group; otherwise $F$ is
orbit-like. Similar results were
obtained by V. Zinde
(see \cite{Zinde74,Zinde77a,Zinde77b,Zinde77c,Zinde78}) for
$X={\mathbb C}^*$.
\vskip5pt

\noindent In this paper, we complete the story for
all non-hyperbolic Riemann surfaces;
one of our main results is as follows.

\begin{Theorem}\label{Thm: tame and orbit-like maps}
For $n>4$, each holomorphic map
$F\colon {\mathcal C}^n({\mathbb T}^2)
\to {\mathcal C}^n({\mathbb T}^2)$
is either tame or orbit-like. In particular,
any automorphism of ${\mathcal C}^n({\mathbb T}^2)$
is tame.

Moreover, $F$ is tame if and only if
the induced endomorphism $F_*$ of the
fundamental group $\pi_1({\mathcal C}^n({\mathbb T}^2))$
is non-abelian, i.e. its image
$F_*(\pi_1({\mathcal C}^n({\mathbb T}^2)))$ is
a non-abelian group.
\end{Theorem}

\noindent This theorem has two immediate corollaries.
The first shows that non-abelian holomorphic endomorphisms of
${\mathcal C}^n({\mathbb T}^2)$ admit a simple classification
up to a holomorphic homotopy.\footnote{That is,
a homotopy within the space of holomorphic mappings.}

\begin{Corollary}\label{Crl: homotopy classification}
For $n>4$, the set
${\mathcal H}({\mathcal C}^n({\mathbb T}^2),
{\mathcal C}^n({\mathbb T}^2))$
of all holomorphic homotopy classes of non-abelian
holomorphic endomorphisms
of ${\mathcal C}^n({\mathbb T}^2)$ is in a natural
one-to-one correspondence with the set
${\mathcal H}({\mathcal C}^n({\mathbb T}^2),\Aut{\mathbb T}^2)$
of all holomorphic homotopy classes
of holomorphic maps ${\mathcal C}^n({\mathbb T}^2)
\to\Aut{\mathbb T}^2$.
\end{Corollary}

\begin{Corollary}\label{Crl: orbit space}
Let $n>4$ and $G=\Aut{\mathcal C}^n({\mathbb T}^2)$\,.
Then the orbits of the natural $G$-action in
${\mathcal C}^n({\mathbb T}^2)$
coincide with the orbits of the diagonal $(\Aut{\mathbb T}^2)$-action in ${\mathcal C}^n({\mathbb T}^2)$.
\end{Corollary}

\noindent Our second main theorem deals with the torus
braid group $B_n({\mathbb T}^2)
=\pi_1({\mathcal C}^n({\mathbb T}^2))$
and its pure braid subgroup
$P_n({\mathbb T}^2)$, which is the fundamental group
of the {\em ordered} configuration space
$$
{\mathcal E}^n({\mathbb T}^2)=\{q=(q_1,...,q_n)\in
({\mathbb T}^2)^n\,|\ q_i\ne q_j \ \forall\, i\ne j\}\,.
$$

\begin{Theorem}\label{Thm: algebraic lifting for torus}
$a)$ Let $n>4$. Then
the pure braid group
$P_n({\mathbb T}^2)\subset B_n({\mathbb T}^2)$
is invariant under any non-abelian endomorphism
$\varphi$ of the whole braid group $B_n({\mathbb T}^2)$,
that is, $\varphi(P_n({\mathbb T}^2))\subseteq
P_n({\mathbb T}^2)$.

$b)$ Let $n>4$ and $n>m$. Then any homomorphism
$\varphi\colon B_n(\mathbb{T}^2)
\to B_m(\mathbb{T}^2)$ is abelian, i.e.
$\varphi(B_n(\mathbb{T}^2))$ is an abelian group.
\end{Theorem}

\noindent For automorphisms of the classical
Artin braid group $B_n({\mathbb C})$, an analogue of part
$(a)$ was proved by Artin himself in 1947
(see \cite{Art47b}).
V. Lin \cite{Lin72c,Lin79,Lin96b,Lin04b,LinSphere}
generalized this result of Artin
to non-abelian endomorphisms of $B_n({\mathbb C})$ and
$B_n({\mathbb{CP}}^1)$, respectively;
the case of $B_n({\mathbb C}^*)$ was handled by V. Zinde
\cite{Zinde77c,Zinde78}. In 1992, N. Ivanov \cite{Ivanov92}
proved a similar result for {\em automorphisms}
of the braid group $B_n(X)$ of any finite type Riemann surface
but for $X={\mathbb{CP}^1}$.
Theorem \ref{Thm: algebraic lifting for torus}$(a)$
completes the story for non-hyperbolic curves.
Analogues of statement $(b)$ were known for
the braid groups $B_n(\mathbb C)$, $B_n({\mathbb{CP}^1})$
and $B_n({\mathbb C}^*)$ (see papers
by V. Lin and V. Zinde quoted above)


\subsection{Plan of the proof of
Theorem \ref{Thm: tame and orbit-like maps}}
\label{ss: 1.2}
First, due to Theorem
\ref{Thm: algebraic lifting for torus}$(a)$,
every continuous non-abelian self-map $F$
fits into a commutative
diagram
\begin{equation}\label{CD: lifting diagram}
\CD
{{\mathcal E}^n(\mathbb{T}^2)} @ > {f} >>
{{\mathcal E}^n(\mathbb{T}^2)}\\
@V{p}VV @VV{p}V\\
{{\mathcal C}^n(\mathbb{T}^2)} @ >> {F}>
{{\mathcal C}^n(\mathbb{T}^2)}\,,
\endCD
\end{equation}
where $p$ is the natural projection
\begin{equation}\label{eq: projection p}
{\mathcal E}^n(\mathbb{T}^2)\ni
q=(q_1,...,q_n)\mapsto \{q_1,...,q_n\}=Q
\in{\mathcal C}^n(\mathbb{T}^2)\,.
\end{equation}
The map $f$ is {\em strictly equivariant} with respect to
the standard action of the symmetric group ${\mathbf S}(n)$
in ${\mathcal E}^n(\mathbb{T}^2)$, meaning that there is
an automorphism $\alpha$ of ${\mathbf S}(n)$ such that
\begin{equation*}
F(\sigma q)=\alpha(\sigma)F(q) \ \
\text{for all} \ q\in{\mathcal E}^n(\mathbb{T}^2) \
\text{and} \ \sigma\in{\mathbf S}(n)\,.
\end{equation*}
Moreover, $f$ is non-constant and
holomorphic whenever $F$ is so.
\vskip5pt

\noindent A torus $\mathbb{T}^2
=\mathbb{C}/\{\rm{a \ lattice \ of \ rank} \ 2\}$
which we deal with is an additive complex Lie group.
To study non-constant strictly equivariant
holomorphic endomorphisms $f$ of the space
${\mathcal E}^n={\mathcal E}^n(\mathbb{T}^2)$,
we start with an explicit description
of all non-constant holomorphic maps
$\lambda\colon{\mathcal E}^n
\to \mathbb{T}^2\setminus \{0\}$.
The set $L$ of all such maps is finite
and separates points of a certain
submanifold $M\subset{\mathcal E}^n$ of complex
codimension $1$. An endomorphism $f$ induces a self-map
$f^*$ of $L$ via $f^*\colon L\ni \lambda
\mapsto f^*\lambda=\lambda\circ f\in L$.

The map $f^*$ carries important information about
$f$. In order to investigate behaviour of $f^*$
and then recover 
this information,
we endow $L$ with the following simplicial structure:
a subset $\Delta^s=\{\lambda_0,...,\lambda_s\}\subseteq L$
is said to be an {\em $s$-simplex} whenever
$\lambda_i-\lambda_j\in L$ for all distinct $i,j$.
The action of ${\mathbf S}(n)$ in ${\mathcal E}^n$
induces a simplicial ${\mathbf S}(n)$-action in
the complex $L$; the orbits of this action may be
exhibited explicitly. On the other hand,
the map $f^*\colon L\to L$ defined above is simplicial and
preserves dimension of simplices; since $f$ is equivariant,
$f^*$ is nicely related to the ${\mathbf S}(n)$ action on $L$.

Studying all these things together, we eventually
come to the conclusion that $f$ is {\em tame},
meaning that there exists an ${\mathbf S}(n)$ invariant
holomorphic map
\begin{equation}\label{eq: invariant map to Aut T2}
t\colon {\mathcal E}^n(\mathbb{T}^2)\to\Aut\mathbb{T}^2
\end{equation}
and a permutation $\sigma\in{\mathbf S}(n)$
such that
\begin{equation}\label{eq: tame map E to E}
f(q)=\sigma t(q)q
=(t(q)q_{\sigma^{-1}(1)},...,t(q)q_{\sigma^{-1}(n)})
\ \ \forall \,
q=(q_1,...,q_n)\in{\mathcal E}^n(\mathbb{T}^2)\,.
\end{equation}
The latter formula implies that the original endomorphism
$F$ is tame.
\vskip7pt

\noindent For the complete proofs, see Section
\ref{Sec: Some algebraic properties of torus braid groups}
(the proof of Theorem \ref{Thm: algebraic lifting for torus}
and some related \linebreak[3] algebraic results), Section
\ref{Sec: Some analytic properties of ordered configuration space}
(a simplicial complex of holomorphic maps 
${\mathcal E}^n(\mathbb{T}^2)\to\mathbb{T}^2\setminus\{0\}$)
and Section
\ref{Sec: Holomorphic mappings of configuration spaces}
(the proof of Theorem \ref{Thm: tame and orbit-like maps}).

\subsection{Notation and definitions}
\label{ss: Notation and definitions}

\noindent For the reader's convenience, we collected here
the main notation and definitions used throughout the paper.
\vskip5pt

\begin{Definition}\label{Def: configuration spaces}
The {\em ordered} and {\em non-ordered
configuration spaces} of a torus $\mathbb{T}^2$
are defined as follows:
$$
\aligned
{\mathcal E}^n(\mathbb{T}^2)
&=\{(q_1,...,q_n)\in (\mathbb{T}^2)^n\,|\
q_i\ne q_j \ \forall\,i\ne j\}\,,\\
{\mathcal C}^n(\mathbb{T}^2)
&=\{Q\subset \mathbb{T}^2\,|\ \#Q=n\}\,.
\endaligned
$$
\noindent The projection
$$
p\colon{\mathcal E}^n(\mathbb{T}^2)\to
{\mathcal C}^n(\mathbb{T}^2)\,, \ \
p(q)=p(q_1,...,q_n)=\{q_1,...,q_n\}=Q\,,
$$
is an ${\mathbf S}(n)$ Galois covering
and we have the exact sequence of the corresponding
fundamental groups
\begin{equation}\label{eq: exact sequence}
1\to\pi_1({\mathcal E}^n(\mathbb{T}^2))
\overset{p_*}{\longrightarrow}
\pi_1({\mathcal C}^n(\mathbb{T}^2))
\to{\mathbf S}(n)\to 1\,.
\end{equation}
The fundamental group $\pi_1({\mathcal C}^n(\mathbb{T}^2))$
is called the {\em torus braid group} and is denoted
by $B_n(\mathbb{T}^2)$. Its normal subgroup
$P_n(\mathbb{T}^2)\Def\pi_1({\mathcal E}^n(\mathbb{T}^2))$
is called the {\em pure torus braid group}.
\end{Definition}

\begin{Definition}\label{Def: abelian homomorphism or map}
A group homomorphism $\varphi\colon G\to H$ is said to be
{\em abelian} if its image $\varphi(G)$ is an abelian
subgroup of $H$. A continuous map $F\colon X\to Y$
of arcwise-connected spaces is called {\em abelian}
if the induced homomorphism $F_*\colon\pi_1(X)\to
\pi_1(Y)$ of the fundamental groups is abelian.
\end{Definition}

\noindent For a complex space $X$, we denote by $\Aut X$
the group of all biholomorphic automorphisms of $X$.
For algebraic $X$, $\Aut_{reg} X$ stands for the group of all
biregular automorphisms.

The group $\Aut\mathbb{T}^2=\Aut_{reg}\mathbb{T}^2$
is a compact complex Lie group isomorphic to a semidirect product
$\mathbb{T}^2\leftthreetimes{\mathbb Z}_k$; here
$k=2$ if the torus $\mathbb{T}^2$ has no complex multiplications;
otherwise, either $k=4$ or $k=6$.

\begin{Definition}\label{Def: tame map}
A holomorphic endomorphism $F$ of
${\mathcal C}^n(\mathbb{T}^2)$ is said to be {\em tame}
if there is a holomorphic map
$T\colon{\mathcal C}^n(\mathbb{T}^2)
\to\Aut\mathbb{T}^2$ such that $F(Q)=T(Q)Q$ for all
$Q\in{\mathcal C}^n(\mathbb{T}^2)$.
\end{Definition}

\begin{Definition}{\label{Def:  orbit-like maps}}
A holomorphic map $F\colon{\mathcal C}^n(\mathbb{T}^2)
\to{\mathcal C}^n(\mathbb{T}^2)$
is said to be {\em  orbit-like} if its
image $F({\mathcal C}^n(\mathbb{T}^2))$ is
contained in one orbit of the diagonal $\Aut \mathbb{T}^2$ action
in ${\mathcal C}^n(\mathbb{T}^2))$.
\end{Definition}

\section{Some algebraic properties of torus braid groups}
\label{Sec: Some algebraic properties of torus braid groups}

\noindent The main aim of this section is to prove
Theorem \ref{Thm: algebraic lifting for torus}.
To this end, we first need to prove some auxiliary results,
which also seem to be of independent interest.

\subsection{Zariski presentation and homomorphism
$B_n\to B_n({\mathbb T}^2)$}
\label{ss: Zariski presentation}
We will use the following presentation
of the torus braid group $B_n({\mathbb T}^2)$ found
by O. Zariski \cite{Zariski37}\footnote{Another presentation
was found by J. Birman \cite{Birman69}.}.
\vskip5pt

{\em Generators:}
\begin{equation}\label{eq: generators of Bn(T2)}
\sigma_1,\dots,\sigma_{n-1}\,; \ a_1,a_2\,;
\end{equation}

{\em relations:}

\begin{align}
&\sigma_i \sigma_j =\sigma_j \sigma_i \hskip65pt {\rm for} \
                   \modo{i-j} \ge 2\,, \    i,j=1,\dots,n-3\,;
                                           \label{rel:Zariski 1}\\
&\sigma_i \sigma_{i+1} \sigma_i  =
       \sigma_{i+1} \sigma_i \sigma_{i+1}  \hskip18pt {\rm for} \
                                        i=1,\dots,n-2\,;\label{rel:Zariski 2}\\
&\sigma_i a_k = a_k \sigma_i \hskip65pt {\rm for} \
                                 k=1,2 \ \ {\rm and} \ \  i=2,\dots,n-1\,;\label{rel:Zariski 3}\\
&(\sigma_1^{-1} a_k)^2 = (a_k \sigma_1^{-1})^2 \hskip25pt {\rm for} \
                                                    k=1,2\,;\label{rel:Zariski 4}\\
&\sigma_1 \dots \sigma_{n-2} \sigma_{n-1}^2 \sigma_{n-2}
                 \dots \sigma_1 = a_1 a_2^{-1} a_1^{-1}
                  a_2\,;\label{rel:Zariski 5}\\
&a_2  \sigma_1^{-1} a_1^{-1} \sigma_1 a_2^{-1} \sigma_1^{-1} a_1
                         \sigma_1 = \sigma_1^2\,.
                         \label{rel:Zariski 6}
\end{align}

\subsection{Normal series in the torus pure braid groups $P_n(\mathbb{T}^2)$}
\label{ss: Normal series in Pn(T2)}
The exact homotopy sequence \eqref{eq: exact sequence} of
the covering $p\colon{\mathcal E}^n(\mathbb{T}^2)\to
{\mathcal C}^n(\mathbb{T}^2)$ may be written as
\begin{equation}\label{eq: algebraic exact sequence}
1\to P_n(\mathbb{T}^2)\overset{p_*}{\longrightarrow}
B_n(\mathbb{T}^2)\overset{\mu}{\longrightarrow}{\mathbf S}(n)\to 1\,,
\end{equation}
where the epimorphism ${\mu}$ is defined by
\begin{equation}\label{eq: definition of mu}
\mu(\sigma_i)=(i,i+1)\ \ {\rm for} \ i=1,...,n-1 \ \ {\rm and} \ \
\mu(a_1)=\mu(a_2)=1\,.
\end{equation}
Let us take a point $(c_1,...,c_n)\in{\mathcal E}^n(\mathbb{T}^2)$
and for each $m=1,...,n-1$ consider the ordered configuration
space ${\mathcal E}^{n-m}(\mathbb{T}^2\setminus\{c_1,...,c_m\})$
of the torus $\mathbb{T}^2$ punctured at $m$ points $c_1,...,c_m$;
a point $q\in{\mathcal E}^{n-m}(\mathbb{T}^2\setminus\{c_1,...,c_m\})$
has $n-m$ pairwise distinct components
$q_{m+1},...,q_n\in\mathbb{T}^2\setminus\{c_1,...,c_m\}$.
Each of these configuration spaces is an Eilenberg--MacLane $K(\pi,1)$
space for its fundamental group $\pi$
(see E. Fadell and L. Neuwirth \cite{FadNeu62}, Corollary 2.2), which is the pure braid group
(on $n-m$ strands) $P_{n-m;m}$ of the punctured torus
$\mathbb{T}^2\setminus\{c_1,...,c_m\}$.
The projections
\begin{align}
&t_1\colon {\mathcal E}^n(\mathbb{T}^2)
\to\mathbb{T}^2\,,\ \ (q_1,...,q_n)\mapsto q_1\,,\label{projection 1}\\
&{\aligned t_{m+1}\colon
{\mathcal E}^{n-m}(\mathbb{T}^2&\setminus\{c_1,...,c_m\})
\to\mathbb{T}^2\setminus\{c_1,...,c_m\}\,, \\
                 &\hskip60pt (q_{m+1},...,q_n)\mapsto q_{m+1}\,, \ \
                                                    1\le m\le n-2\,,
\endaligned}
\label{projection m}
\end{align}
define locally trivial, smooth fibrings\footnote{See \cite{FadNeu62}
for the case of arbitrary Riemann surfaces.} with fibres isomorphic to
${\mathcal E}^{n-1}(\mathbb{T}^2\setminus\{c_1\})$ and
${\mathcal E}^{n-m-1}(\mathbb{T}^2\setminus\{c_1,...,c_{m+1}\})$,
respectively. The final segments of the exact homotopy sequences
of these fibrings look, respectively, as
$$
1\to P_{n-1;1}\to P_n(\mathbb{T}^2)\to{\mathbb Z}^2\to 1
$$
and
$$
1\to P_{n-m-1;m+1}\to P_{n-m;m}\to{\mathbb F}_m\to 1\,,
$$
where ${\mathbb F}_m$ stands for a free group of rank $m$.
This leads immediately to the following statement,
which is an analogue of a well-known Markov theorem for Artin
pure braid groups, proved in \cite{Mar45}.

\begin{Proposition}\label{Prp: normal series}
The subgroups $P_{n-s;s}$ fit into the normal series:
$$
\{1\}\subset P_{1;n-1}\subset\cdots\subset P_{n-m-1;m+1}\subset
P_{n-m;m}\subset\cdots\subset P_{n-1;1}
\subset P_{n;0} = P_n(\mathbb{T}^2)
$$
with the factors
\begin{equation}\label{eq: factors}
\aligned
P_{1;n-1}/\{1\}&=P_{1;n-1}=\pi_1(\mathbb{T}^2\setminus\{a_1,...,a_{n-1}\})
\cong {\mathbb F}_{n-1}\,, \ ...\,, \\
&P_{n-m;m}/P_{n-m-1;m+1}\cong {\mathbb F}_m\,, \ ...\,, \\
&P_{n-1;1}/P_{n-2;2}\cong {\mathbb F}_2\,, \ \ \
P_n(\mathbb{T}^2)/P_{n-1;1}\cong {\mathbb Z}^2\,.
\endaligned
\end{equation}
In particular, any non-trivial subgroup
$H\subseteq P_n(\mathbb{T}^2)$
admits a non-trivial homomorphism to ${\mathbb Z}$.
\end{Proposition}


\noindent The last sentence of the proposition implies:

\begin{Proposition}\label{Prp: homomorphisms to Pn(T2)}
A group $G$ with finite abelianization $G/G'=G/[G,G]$
has no non-trivial homomorphisms to the pure torus braid group
$P_n(\mathbb{T}^2)$.
\end{Proposition}

\noindent Recall that the Artin braid group
$B_n=\pi_1({\mathcal C}^n({\mathbb C}))$
has generators $\sigma_1,...,\sigma_{n-1}$
and the system of defining relations \eqref{rel:Zariski 1},
\eqref{rel:Zariski 2}. Thus, we have a natural homomorphism
$i\colon B_n\to B_n({\mathbb T}^2)$ sending all generators
$\sigma$'s of $B_n$
to the elements of the same name in $B_n({\mathbb T}^2)$ (according to
N. Ivanov \cite{Ivanov92}, $i$ is injective; however,
we will not use this fact).

\begin{Lemma}\label{Lm: composition of homomorphisms}
Let $n>4$ and let $\varphi\colon B_n({\mathbb T}^2)\to B_m({\mathbb T}^2)$
be a homomorphism such that the composition
\begin{equation}\label{eq: composition mu phi i}
\varPhi=\mu\circ\varphi\circ i
\colon B_n\stackrel{i}{\longrightarrow}
B_n(\mathbb{T}^2)\stackrel{\varphi}{\longrightarrow}
B_m(\mathbb{T}^2) \stackrel{\mu}{\longrightarrow}{\mathbf S}(m)
\end{equation}
is an abelian homomorphism. Then $\varphi$ is abelian.
In particular, $\varphi$ is abelian. whenever $\mu\circ\varphi$ is so. 
\end{Lemma}

\begin{proof}
Let $\varPhi'\colon B_n'\to {\mathbf S}(m)$ denote the restriction of
$\varPhi$ to the commutator subgroup $B_n'$ of the group $B_n$.
Since $\varPhi$ is abelian, $\varPhi'$ is trivial
and hence $\varphi(i(B_n'))\subseteq \Ker\mu= P_m(\mathbb{T}^2)$.
By Gorin--Lin Theorem \cite{GorLin69}, for $n>4$ the group $B_n'$
is perfect\footnote{That is, $B_n'$ coincides with its commutator subgroup
$B_n''=[B_n',B_n']$.}, and Proposition 
\ref{Prp: homomorphisms to Pn(T2)}
shows that $\varphi(i(B_n'))=1$. Hence $\varphi\circ i$ is abelian,
and relations $\sigma_k\sigma_{k+1}\sigma_k=\sigma_{k+1}\sigma_{k}\sigma_{k+1}$
imply $\varphi(i(\sigma_1))=...=\varphi(i(\sigma_{n-1}))$.
By definition, $i$ sends the generators $\sigma_1,...,\sigma_{n-1}$
of $B_n$ into the corresponding generators $\sigma_1,...,\sigma_{n-1}$
of $B_n(\mathbb{T}^2)$; thus, the latter ones satisfy
\begin{equation}\label{eq: all phi(sigma) coincide}
\varphi(\sigma_1)=...=\varphi(\sigma_{n-1})\,.
\end{equation}
Together with
\eqref{rel:Zariski 3} this shows that
\begin{equation}\label{eq: phi(ak) and phi(sigmaj) commute}
\varphi(a_k)\varphi(\sigma_j)=\varphi(\sigma_j)\varphi(a_k) \ \
{\rm for \ all} \ j,k\,.
\end{equation}
By \eqref{eq: all phi(sigma) coincide} and \eqref{rel:Zariski 5},
\begin{equation}\label{eq: aux equation1}
(\varphi(\sigma_1))^{2(n-1)}=
\varphi(a_1) \varphi(a_2)^{-1} \varphi(a_1)^{-1} \varphi(a_2)\,,
\end{equation}
and by \eqref{eq: phi(ak) and phi(sigmaj) commute} and \eqref{rel:Zariski 6},
\begin{equation}\label{eq: aux equation2}
\varphi(a_2) \varphi(a_1)^{-1} \varphi(a_2)^{-1} \varphi(a_1)
= \varphi(\sigma_1)^2\,.
\end{equation}
Relations \eqref{eq: phi(ak) and phi(sigmaj) commute} and
\eqref{eq: aux equation2} imply
\begin{equation}\label{eq: aux equation3}
\varphi(\sigma_1)^2 =
\varphi(a_2)^{-1} \varphi(a_1)\varphi(a_2)\varphi(a_1)^{-1}  \,.
\end{equation}
Multiplying \eqref{eq: aux equation1}
and \eqref{eq: aux equation3} we obtain
\begin{equation*}
1=  (\varphi (\sigma_1))^{2(n-1)}
(\varphi (\sigma_1))^{2}=(\varphi (\sigma_1))^{2n}\,.
\end{equation*}
However, by E. Fadell and L. Neuwirth \cite{FadNeu62},
Theorem $8$, the torus braid group
$B_m(\mathbb{T}^2)$ has no elements of a finite
order, this implies $\varphi (\sigma_1)=1$. Thus
\begin{equation*}
\varphi(a_2)^{-1} \varphi(a_1)\varphi(a_2)\varphi(a_1)^{-1}=1,
\end{equation*}
and $\varphi$ is an abelian homomorphism.
\end{proof}

\subsection{Proof of Theorem
                {\ref{Thm: algebraic lifting for torus}}}
\label{ss: Proof of Theorem about algebraic lifting}
Let $n>4$ and let $\varphi$ be a non-abelian endomorphism
of $B_n({\mathbb T}^2)$. By Lemma \ref{Lm: composition of homomorphisms},
the homomorphism $\varPhi=\mu\circ\varphi\circ i\colon B_n
\to{\mathbf S}(n)$ in \eqref{eq: composition mu phi i} is
non-abelian.
Thus, by V. Lin's theorem (see \cite{Lin79}, Section 4;
a complete proof may be found in \cite{Ivanov92,Lin96b,Lin04b}),
$\varPhi$ coincides with the standard epimorphism
$B_n\to{\mathbf S}(n)$
up to an automorphism of ${\mathbf S}(n)$. It follows
that the homomorphism $\mu\circ\varphi
\colon B_n({\mathbb T}^2)\to{\mathbf S}(n)$ is surjective.
N. Ivanov (see \cite{Ivanov92}, Theorem 1)
proved that for $n>4$, any non-abelian
homomorphism $B_n({\mathbb T}^2)\to{\mathbf S}(n)$
whose image is a primitive permutation group on $n$ letters
coincides with the standard epimorphism
$\mu\colon B_n({\mathbb T}^2)\to{\mathbf S}(n)$
up to an automorphism of ${\mathbf S}(n)$. Thus,
$\Ker(\mu\circ\varphi)=P_n({\mathbb T}^2)=\Ker\mu$, which
implies $\varphi(P_n({\mathbb T}^2))\subseteq
\Ker\mu=P_n({\mathbb T}^2)$ thus proving part $(a)$ of the
theorem.\footnote{In fact, we proved a slightly stronger
property $\varphi^{-1}(P_n({\mathbb T}^2))=P_n({\mathbb T}^2)$.}

To prove part $(b)$, we use another Lin's theorem (see quotations above), which says that
for $n>\max(m,4)$ any homomorphism $B_n\to {\mathbf S}(m)$ is abelian; it follows that the homomorphism
$\varPhi=\mu\circ\varphi\circ i\colon B_n\to{\mathbf S}(m)$
in \eqref{eq: composition mu phi i} is abelian.
By Lemma \ref{Lm: composition of homomorphisms}, $\varphi$ is abelian.
\hfill $\square$
\vskip7pt

\noindent We skip the proofs of the next two
results, which will not be used in what follows.

\begin{Theorem}\label{lm: torus commutator subgroup}
For $n \ge 5$, the abelianization 
$B'_n(\mathbb{T}^2)/[B'_n(\mathbb{T}^2),B'_n(\mathbb{T}^2)]$
of the commutator subgroup $B'_n(\mathbb{T}^2)$
of the torus braid group $B_n(\mathbb{T}^2)$ is isomorphic to
$\mathbb{Z}_n=\mathbb{Z}/n\mathbb{Z}$.
\end{Theorem}

\begin{Theorem}
\label{cor: algebraic lifting for torus small dimension}
Let $3\le m \le n \le 4$ and let 
$\varphi\colon B_n(\mathbb{T}^2)\to B_m(\mathbb{T}^2)$ be a
non-abelian homomorphism such that $\mu\circ\varphi$ is non-abelian.
Then $\varphi(P_n(\mathbb{T}^2)) \subseteq P_m(\mathbb{T}^2)$.
\end{Theorem}


\section{Ordered configuration spaces}
\label{Sec: Some analytic properties of ordered configuration space}
\noindent Here we establish some analytic properties of
ordered configuration spaces.

\subsection{A Cartesian product structure
in ${\mathcal E}^n(\mathbb{T}^2)$}
\label{ss: A cartesian product structure in En}
The torus $\mathbb{T}^2$ acts in ${\mathcal E}^n(\mathbb{T}^2)$
by translations: $q=(q_1,...,q_n)\mapsto (q_1+t,...,q_n+t)$,
$t\in\mathbb{T}^2$. Any orbit of this action is isomorphic
to $\mathbb{T}^2$ and intersects the submanifold
$$
\tilde{M}_0=\{q=(q_1,...,q_{n-1},0)
\in{\mathcal E}^n(\mathbb{T}^2)\}
\cong {\mathcal E}^{n-1}(\mathbb{T}^2\setminus \{0\})
$$
at a single point. This gives the Cartesian decomposition
${\mathcal E}^n(\mathbb{T}^2)=\mathbb{T}^2\times \tilde{M}_0$, defined by the following maps
\begin{equation*}
\aligned
{\mathcal E}^n(\mathbb{T}^2)\ni(q_1,...,q_n)\mapsto
(q_n,(q_1-q_n,...,q_{n-1}-q_n,0))\in\mathbb{T}^2\times \tilde{M}_0\,,\\
\mathbb{T}^2\times \tilde{M}_0\ni(t,(q_1,...,q_{n-1},0))\mapsto
(q_1+t,...,q_{n-1}+t,t)\in{\mathcal E}^n(\mathbb{T}^2)\,.
\endaligned
\end{equation*}
It is easily seen that $\tilde{M}_0$ is a non-singular
irreducible affine algebraic variety; in particular, it is
a Stein manifold, whereas ${\mathcal E}^n(\mathbb{T}^2)$
is not so.

\subsection{Holomorphic mappings
${\mathcal E}^n(\mathbb{T}^2)\to\mathbb{T}^2\setminus\{0\}$}
\label{ss: Holomorphic mappings En to T2-0}
We continue to execute our plan sketched in Section
\ref{ss: 1.2} with the following

\begin{Theorem}\label{Thm: description of all differences}
Any non-constant holomorphic map $
f:\mathcal{E}^n(\mathbb{T}^2)
\rightarrow\mathbb{T}^2\setminus\{0\}$
is of the form
\begin{equation*}
f(q_1,q_2,\dots,q_n)= \mathfrak{m}(q_i-q_j) \
{\rm with \ some} \ i\ne j\,,
\end{equation*}
where $\mathfrak{m}\colon
\mathbb{T}^2\to \mathbb{T}^2$ is either the
identity or a complex multiplication.\footnote{Of course,
the latter may happen only if our torus admits such a
multiplication
(by a complex multiplication we mean an automorphism
of the torus $\mathbb{T}^2$ induced by multiplication on the
complex line by a complex number $\ne 1,-1$).}
\end{Theorem}

\begin{Notation}\label{Not: complex multiplications}
We denote by $\mathfrak{M}$ the finite cyclic group
of automorphisms of $\mathbb{T}^2$ consisting of $\pm$identity
and all complex multiplications in $\mathbb{T}^2$ (if they exist).
$\mathfrak{M}$ is isomorphic either to ${\mathbb Z}_2$ or to
${\mathbb Z}_4$ or to ${\mathbb Z}_6$. Let $\mathfrak{M}_+$
consists of all $\mathfrak{m}\in\mathfrak{M}$ with
$0 \le \Arg \mathfrak{m} <\pi$ (i.e. $\mathfrak{M}_+$ contains
$1$, $2$ or $3$ elements).
\end{Notation}

\noindent To prove the theorem we need some preparation.
The proof of the following well-known statement is due to H. Huber
\cite{Huber53}, \S 6, Satz 2; it may also be
found in Sh. Kobayashi's book \cite{Kob70}, Chapter VI,
Sec. 2, remarks after Corollary 2.6.
\vskip5pt

\noindent{\caps Claim A.} {\sl Let $M$ and $N$ be compact Riemann
surfaces and $A\subset M$ and $B\subset N$ be finite sets.
Suppose that the domain $N\setminus B$ is hyperbolic, meaning that
its universal covering is isomorphic to the unit disc. Then any
holomorphic map $M\setminus A\to N\setminus B$
extends to a holomorphic map $M\to N$.}
\vskip5pt

\noindent The following facts are also very well known.
\vskip5pt

\noindent{\caps Claim B.} {\sl Any non-constant holomorphic self-map
$\mu\colon{\mathbb T}^2\to{\mathbb T}^2$ is regular,
has no critical points and hence is an unbranched covering
of some finite degree $k$.
If $\mu(0)=0$ then $\mu$ is also a regular endomorphism of the
compact algebraic group ${\mathbb T}^2$ and its kernel $\Ker\mu$
is a subgroup of order $k$. 
Every holomorphic self-map of ${\mathbb T}^2$ is of the form
$q\mapsto\sum\limits_{{\mathfrak m}
\in{\mathfrak M}}k_{\mathfrak m}{\mathfrak m}(q-a)$, and any
automorphism is of the form
$q\mapsto{\mathfrak m}(q-a)$, where ${\mathfrak m}\in{\mathfrak M}$,
$k_{\mathfrak m}\in{\mathbb Z}$, and $a\in{\mathbb T}^2$.}
\vskip5pt

\noindent {\em Explanation.} The graph
$G_\mu=\{(p,q)\in{\mathbb T}^2\times{\mathbb T}^2\,|\ q=\mu(p)\}$
of $\mu$ is an analytic subset of the projective variety
${\mathbb T}^2\times{\mathbb T}^2$. By general Chow's theorem
(see \cite{Chow} or \cite{GR}, Chapter V, Theorem D7), $G_\mu$
is a (non-singular) projective variety. Of course, the projections
$\alpha\colon G_\mu\ni(p,q)\mapsto p\in{\mathbb T}^2$ and
$\beta\colon G_\mu\ni(p,q)\mapsto q\in{\mathbb T}^2$ are regular.
Hence, the inverse map $\alpha^{-1}\colon{\mathbb T}^2\ni p
\mapsto (p,\mu(p))\in G_\mu$ is regular, and the composition
$\mu=\beta\circ\alpha^{-1}$ is regular, as well.

To see that $\mu$ has no critical points, suppose to the contrary
that the set of all critical values $C$ contains $m\ge 1$ points.
The restriction of $\mu$ to the complement of $\mu^{-1}(C)$ defines
an unbranched covering ${\mathbb T}^2\setminus\mu^{-1}(C)
\to{\mathbb T}^2\setminus C$ of some finite degree $k$. The
Riemann--Hurwitz formula for Euler characteristics states that
$\chi({\mathbb T}^2\setminus\mu^{-1}(C))
=k\chi({\mathbb T}^2\setminus C)$, i.e. $\#(\mu^{-1}(C))=km$.
Since $\#(\mu^{-1}(c))\le k-1$ for all $c\in C$, we have
$km=\#(\mu^{-1}(C))\le (k-1)m$, which is impossible.

Assuming that $\mu(0)=0$, define a map
$\nu\colon{\mathbb T}^2\times{\mathbb T}^2\to{\mathbb T}^2$
by
$$
\nu(p,q)=\mu(p+q)-\mu(p)-\mu(q)\,.
$$
Any loop in ${\mathbb T}^2\times{\mathbb T}^2$ can be deformed
to the subset $({\mathbb T}^2\times\{0\})\cup(\{0\}\times{\mathbb T}^2)$,
where $\nu=0$. Thus, $\nu$ induces the trivial homomorphism
of the fundamental groups and lifts to a holomorphic map
$\tilde\nu\colon{\mathbb T}^2\times{\mathbb T}^2\to{\mathbb C}$,
which, by the maximum principle, must be constant. It follows that
$\nu=\const=0$ and $\mu$ is a holomorphic group endomorphism
of $\mathbb{T}^2$. (Another argument may be found
in \cite{Clemens80}, Chapter 3, Section 3.1.) 

The general form of automorphisms and endomorphisms of ${\mathbb T}^2$
may be found in \cite{Farkash}, Chapter V, Section V.4.7, and
\cite{Lang73}, Chapter 1, Section 5, respectively.

\hfill $\square$
\vskip5pt

\begin{Lemma}\label{Lm: holomorphic maps T2-a to T2-b}
For any $a,b\in\mathbb{T}^2$ every non-constant holomorphic
map $\lambda\colon\mathbb{T}^2\setminus\{a\}\to\mathbb{T}^2\setminus\{b\}$
extends to a biregular automorphism of $\mathbb{T}^2$ sending
$a$ to $b$.
\end{Lemma}

\begin{proof}
By Claims A and B, $\lambda$ extends to an unbranched holomorphic covering map
$\tilde\lambda\colon\mathbb{T}^2\to\mathbb{T}^2$ of some finite degree $k$.
Clearly $\tilde\lambda^{-1}(a)=\{b\}$; hence, $k=1$ and $\tilde\lambda$
is an automorphism.
\end{proof}

\begin{Definition}\label{Def: exceptional configuration}
A configuration $a=(a_1,...,a_m)\in{\mathcal E}^m({\mathbb T}^2)$
is said to be {\em exceptional} if there exist $i\ne j$
and a non-constant
holomorphic self-map $\lambda\colon{\mathbb T}^2\to{\mathbb T}^2$
such that $\lambda(a_i)=\lambda(a_j)$ and
$\lambda^{-1}(\lambda(a_i))\subseteq\{a_1,...,a_m\}$.
\end{Definition}

\begin{Lemma}\label{Lm: exceptional configuration}
$a)$ The set $A$ of all exceptional configurations
$a=(a_1,...,a_m)\in{\mathcal E}^m({\mathbb T}^2)$ is contained
in a subvariety $M\subset{\mathcal E}^m({\mathbb T}^2)$
of codimension $1$.
\vskip4pt

$b)$ For any non-exceptional configuration
$(a_1,...,a_m)\in{\mathcal E}^m({\mathbb T}^2)$,
every non-cons\-tant holomorphic map
$\lambda\colon\mathbb{T}^2\setminus\{a_1,...,a_m\}
\to\mathbb{T}^2\setminus\{0\}$
extends to a biregular automorphism of $\mathbb{T}^2$, sending
a certain $a_i$ to $0$.
\end{Lemma}

\begin{proof}
$a)$ Let $N$ denote the union of all finite subgroups
of order $\le m$ in $\mathbb{T}^2$; this set is finite.
Set $M=\{(a_1,...,a_m)\in{\mathcal E}^m({\mathbb T}^2)
|\,a_j-a_i\in N \ {\rm for \ some} \ i\ne j\}$; then
$M$ is a subvariety in ${\mathcal E}^m({\mathbb T}^2)$
of codimension $1$. We shall show that $A\subseteq M$.

Let $a=(a_1,...,a_m)\in A$ and let $i,j$, and
$\lambda$ be as in Definition
\ref{Def: exceptional configuration}.
Set $\mu(t)=\lambda (t+a_i)-\lambda(a_i)$, $t\in\mathbb{T}^2$.
Then $\mu(0)=0$ and, by Claim B,
$\mu$ is a group homomorphism with finite kernel
$\#\Ker\mu$. If $t\in\Ker\mu$, then $\lambda(t+a_i)=\lambda(a_i)$,
$t+a_i\in\lambda^{-1}(\lambda(a_i))\subseteq
\{a_1,...,a_m\}$ and $t\in\{a_1-a_i,...,0,...,a_m-a_i\}$;
that is, $\Ker\mu\subseteq\{a_1-a_i,...,0,...,a_m-a_i\}$.
In particular, $\#\Ker\mu\le m$ and hence $\Ker\mu\subseteq N$.
Since $\mu(a_j-a_i)=0$, we have $a_j-a_i\in N$ and $a\in M$.
\vskip5pt

$b)$ Let $a=(a_1,...,a_m)\notin A$.
By Claims A and B, any non-constant
holomorphic map 
$\lambda\colon\mathbb{T}^2\setminus\{a_1,...,a_m\}
\to\mathbb{T}^2\setminus\{0\}$ extends to
a finite unbranched holomorphic covering map
$\tilde\lambda\colon\mathbb{T}^2\to\mathbb{T}^2$.
Clearly $\tilde\lambda^{-1}(0)\subseteq\{a_1,...,a_m\}$;
in particular, $\tilde\lambda(a_i)=0$ for a certain $i$.
Since $a\notin A$, we have
$\tilde\lambda(a_j)\ne 0$ for all $j\ne i$;
this means that $\tilde\lambda^{-1}(0)=\{a_i\}$
and $\deg\tilde\lambda=1$.
\end{proof}

\noindent {\em Proof of Theorem}
\ref{Thm: description of all differences}
First, we notice that any holomorphic map
$\mathbb{T}^2\to\mathbb{T}^2\setminus\{0\}$
is constant, since it lifts to a holomorphic
map ${\mathbb C}\to{\mathbb D}
=\{\zeta\in{\mathbb C}\,|\ |\zeta|<1\}$ of the universal coverings,
which is constant due to Liouville's theorem.
\vskip5pt

\noindent The proof of the theorem is by induction on $n$,
starting with $n=2$. For $a\in\mathbb{T}^2$, denote by
$\lambda_a=\lambda(\cdot,a)$ the restriction of $\lambda$
to the fibre $p^{-1}(a)=\mathbb{T}^2\setminus\{a\}$
of the projection $p\colon\mathcal{E}^2(\mathbb{T}^2)\ni
(q_1,q_2)\mapsto q_2\in\mathbb{T}^2$.

The set $S=\{a\in\mathbb{T}^2\,|\ \lambda_a \ 
\text{\rm is a constant map}\}$ is finite. Because,
otherwise, by the uniqueness theorem,
$\lambda_a=\const$ for all $a\in\mathbb{T}^2$, \
$\lambda=\lambda(q_1,q_2)$ does not depend on $q_1$ and
may be considered as a holomorphic map $\mathbb{T}^2
\to\mathbb{T}^2\setminus\{0\}$, which must be constant;
however, this implies $\lambda=\const$, contradicting our
assumption.

By Lemma \ref{Lm: holomorphic maps T2-a to T2-b},
for any $a\notin S$ the map $\lambda_a\colon
\mathbb{T}^2\setminus\{a\}\to\mathbb{T}^2\setminus\{0\}$
extends to a biholomorphic automorphism of $\mathbb{T}^2$
sending $a$ to $0$. This means that $\lambda_a$
is of the form $\lambda_a(t) = \mathfrak{m}(t-a)$,
with some $\mathfrak{m}=\mathfrak{m}_a\in\mathfrak{M}$.
Thus, for all $(q_1,q_2)$ in the connected, everywhere dense
set $\mathcal{E}^2(\mathbb{T}^2)\setminus p^{-1}(S)$
we have
\begin{equation}\label{eq: almost everywhere}
\lambda(q_1,q_2)=\mathfrak{m}(q_1-q_2)
\end{equation}
with a certain $\mathfrak{m}=\mathfrak{m}_{q_2}\in\mathfrak{M}$.
Since $\mathfrak{M}$ is finite, the element
$\mathfrak{m}=\mathfrak{m}_{q_2}\in\mathfrak{M}$ on
the right hand side of \eqref{eq: almost everywhere} cannot depend
on $q_2$ and the latter formula holds true for the whole of
$\mathcal{E}^2(\mathbb{T}^2)$, which completes the proof for $n=2$.
\vskip5pt

\noindent Assume that the theorem is already proved for some $n=m-1\ge 2$.
For $a=(a_2,...,a_m)\in{\mathcal E}^{m-1}(\mathbb{T}^2)$, denote by
$\lambda_a=\lambda(\cdot,a_2,\dots,a_m)$ the restriction of $\lambda$
to the fibre
$p^{-1}(a)=\mathbb{T}^2\setminus\{a_2,...,a_m\}$
of the projection $p\colon\mathcal{E}^m(\mathbb{T}^2)\ni
(q_1,q_2,...,q_m)\mapsto (q_2,...,q_m)\in\mathcal{E}^{m-1}(\mathbb{T}^2)$.

It is clear that $S\Def\{a\in{\mathcal E}^{m-1}(\mathbb{T}^2)\,
|\ \lambda_a \ \text{\rm is a constant map}\}$ is an analytic subset
of ${\mathcal E}^{m-1}(\mathbb{T}^2)$, and
either $(i)$ $S = {\mathcal E}^{m-1}(\mathbb{T}^2)$
or $(ii)$ $\dim_{\mathbb C}S\le m-2$.

In case $(i)$, $\lambda=\lambda(q_1,...,q_m)$ does not depend on $q_1$ and
may be considered as a holomorphic map ${\mathcal E}^{m-1}(\mathbb{T}^2)
\to\mathbb{T}^2\setminus\{0\}$; by the induction hypothesis,
$\lambda$ is of the desired form.

Let us consider case $(ii)$.
By Lemma \ref{Lm: exceptional configuration}$(a)$,
the set $A$ of all exceptional configurations
is contained in a subvariety $M\subset {\mathcal E}^{m-1}(\mathbb{T}^2)$
of dimension $m-2$.
Let $a\in {\mathcal E}^{m-1}(\mathbb{T}^2)\setminus (S\cup M)$. Then
$\lambda_a\colon \mathbb{T}^2 \setminus \{a_2,...,a_m\}$
is a non-constant map.
By Lemma \ref{Lm: exceptional configuration}$(b)$, $\lambda_a$ extends to
a biholomorphic map $\tilde\lambda_a\colon\mathbb{T}^2\to\mathbb{T}^2$.
Therefore, $\tilde\lambda_a(t) = \mathfrak{m}(t-a_i)$,
with some $\mathfrak{m}=\mathfrak{m}_a\in\mathfrak{M}$ and $i=i_a$.
Thus, for all $(q_1,...,q_m)$ in the connected, everywhere dense
set $\mathcal{E}^m(\mathbb{T}^2)\setminus p^{-1}(S\cup M)$
we have
\begin{equation}\label{eq: almost everywhere gen case}
\lambda(q_1,...,q_m)=\mathfrak{m}(q_1-q_i)
\end{equation}
with certain $\mathfrak{m}=\mathfrak{m}_{q}\in\mathfrak{M}$ and $i=i_q$.
Since $\mathfrak{M}$ is finite,
$\mathfrak{m}$ and $i$ do not depend
on $q$, and \eqref{eq: almost everywhere gen case} holds true on the
whole of
$\mathcal{E}^m(\mathbb{T}^2)$, which completes the step of induction
thus proving the theorem.
\hfill $\square$
\vskip7pt

\begin{Definition}\label{Def: differences}
For any $\mathfrak{m} \in {\mathfrak M}_+$ and 
$i\ne j \in \{1,\dots,n\}$, the map
\begin{equation}\label{eq: differences}
\aligned
&e_{\mathfrak{m};i,j}\colon\mathcal{E}^n(\mathbb{T}^2)\to
\mathbb{T}^2\setminus\{0\}\,,\\
&e_{\mathfrak{m};i,j}(q)=\mathfrak{m}(q_i-q_j) \ {\rm for } \
q=(q_1,...,q_n)\in \mathcal{E}^n(\mathbb{T}^2)\,,
\endaligned
\end{equation}
is called a {\em difference}.
For a difference $\mu=e_{\mathfrak{m};i,j}$, 
the unordered pair of variables $\{q_i,q_j\}$ is
called the {\em support} of $\mu$ and
the automorphism $\mathfrak{m}\in\mathfrak{M}_+$ is called
the {\em marker} of $\mu$.
We denote them by $\supp\mu$ and $\mathfrak{m}_\mu$ respectively.
It follows from Theorem \ref{Thm: description of all differences}
that any non-constant holomorphic map $\mu:\mathcal{E}^n(\mathbb{T}^2)\to
\mathbb{T}^2\setminus\{0\}$ admits a unique representation in the form
of a difference, that is, $\mu=e_{\mathfrak{m};i,j}$ for some
uniquely defined $\mathfrak{m}\in\mathfrak{M}_+$ and $i,j\in\{1,\dots,n\}$. 
\end{Definition}


\subsection{A simplicial structure on the set of differences}
\label{ss: Simplicial structure in L(En(T2))}
The set of all non-constant holomorphic maps of a complex space to a
punctured torus carries a natural simplicial structure\footnote{Compare
to \cite{LinSphere}}.

\begin{Definition}\label{Def: difference}
For a connected complex space $Y$, let $L(Y)$ denote the set of all
non-constant holomorphic maps $\mu\colon Y \to \mathbb{T}^2
\setminus\{0\}$. For $\mu,\nu\in L(Y)$, we say that $\nu$
is a {\em proper reminder} of $\mu$ and write
$\nu\mid\mu$ if $\mu-\nu\in L(Y)$. This relation is symmetric, i.e.
$\nu\mid\mu$ is equivalent to $\mu\mid\nu$.

We define the graph $\Gamma(Y)$ with the set of vertices $L(Y)$
as follows: $\{\mu,\nu\}$ is an {\em edge} connecting $\mu$ and $\nu$
whenever $\mu\mid\nu$. We denote by $L_{\vartriangle}(Y)$ the
{\em flag complex} of the graph $\Gamma(Y)$. This means that a subset
$\Delta^m=\{\mu_0,...,\mu_m\}\subseteq L(Y)$ is said to be
an $m$-{\em simplex} in $L_{\vartriangle}(Y)$ if $\{\mu_i,\mu_j\}$
is an edge in $\Gamma(Y)$ for all $i\ne j$, i.e.
$\mu_i\mid\mu_j$ for all $i\ne j$.
\end{Definition}

\begin{Lemma}\label{Lm: simplicial map f*: L(Y) to L(Z)}
Let $f\colon Z\to Y$ be a holomorphic map of
connected complex spaces.
Suppose that for each $\lambda\in L(Y)$ the composition
\begin{equation}
f^*(\lambda)\Def\lambda\circ f\colon Z
\overset{f}\longrightarrow Y
\overset{\lambda}\longrightarrow{\mathbb T}^2\setminus\{0\}
\end{equation}
is non-constant. Then
\begin{equation}\label{eq: simplicial map f*: L(Y) to L(Z)}
f^*\colon L(Y)\ni\lambda\mapsto\lambda\circ f\in L(Z)
\end{equation}
is a simplicial map whose restriction to any
simplex $\Delta\subseteq L(Y)$ is injective. In particular,
$f^*$ preserves dimension of simplices. 
\end{Lemma}

\begin{proof}
For any $\lambda\in L(Y)$, $f^*(\lambda)$ is a non-constant
holomorphic map to ${\mathbb T}^2\setminus\{0\}$; hence
$f^*(\lambda)\in L(Z)$. If $\mu,\nu\in L(Y)$ and $\mu\mid\nu$,
then $\lambda=\mu-\nu\in L(Y)$
and $f^*(\mu)-f^*(\nu)=f^*(\mu-\nu)=f^*(\lambda)\in L(Z)$;
consequently, $f^*(\mu)\mid f^*(\nu)$. This implies that the
map $f^*$ is simplicial and injective on any simplex. 
\end{proof}

\begin{Remark}\label{Rmk: dominant map}
Condition $f^*(\lambda)\ne\const$ for all $\lambda\in L(Y)$
is certainly fulfilled for any regular dominant map $f\colon Y\to Z$
of non-singular irreducible algebraic varieties.
\end{Remark}

\noindent We need to study some properties of complexes
$L_{\vartriangle}(Y)$ for the case when $Y$ is
an ordered configuration space of a torus.
Notice that according to Theorem
\ref{Thm: description of all differences}, the set
$L(\mathcal{E}^n(\mathbb{T}^2))$ coincides with the set of all
differences on $\mathcal{E}^n(\mathbb{T}^2)$. \vskip7pt


\begin{Lemma}\label{lm: special property of the L complex}
Let $\{\mu_0,...,\mu_s\}\in
L_{\vartriangle}({\mathcal E}^n(\mathbb{T}^2))$
be an $s$-simplex. Then
$\mathfrak{m}_{\mu_{i}} =\mathfrak{m}_{\mu_{j}}$, \
$\card(\supp\mu_{i} \cap \supp\mu_{j}) = 1$
for all $i\ne j$, and $\card(\supp\mu_0\cap\cdots\cap\supp\mu_s)= 1$.
\end{Lemma}
\begin{proof}
Let $i\ne j$ and let
$\mu_{i}=\mathfrak{m}_i(q_{i'}-q_{i''})$ and
$\mu_{j}=\mathfrak{m}_j(q_{j'}-q_{j''})$.
Since $\mu_i\mid\mu_j$, we must have
$\mu_i-\mu_j=\mathfrak{m}(q_{k'}-q_{k''})$
for some $\mathfrak{m}\in\mathfrak{M}_+$ and $k'\ne k''$. Thus,
$\mathfrak{m}_i(q_{i'}-q_{i''})-\mathfrak{m}_j(q_{j'}-q_{j''})=
\mathfrak{m}(q_{k'}-q_{k''})$.
The latter relation can be fulfilled only if
either $\mathfrak{m}_i q_{i'}-\mathfrak{m}_j q_{j'}=0$ or
$\mathfrak{m}_i q_{i''}-\mathfrak{m}_j q_{j''}=0$.
This implies $\mathfrak{m}_i = \mathfrak{m}_j$ and
we have
\begin{equation}\label{eqvs in proper reminder case}
{\rm either} \hskip20pt i'=j' \hskip20pt {\rm or} \hskip20pt i''=j''.
\end{equation}
If $s=1$ we have finished the proof. If $s>2$, then the property
$\#(\supp\mu_i\cap\supp\mu_j)=1$ implies immediately
that $\#(\supp\mu_0\cap\cdots\cap\supp\mu_s)=1$.
For $s=2$ we have
$$
\mu_0 = \mathfrak{m}(q_{i'}-q_{i''})\,, \ \
\mu_1 = \mathfrak{m}(q_{j'}-q_{j''})\,, \ \
\mu_2 = \mathfrak{m}(q_{k'}-q_{k''})\,.
$$
Since $\mu_0 \mid \mu_1$, $\mu_1 \mid \mu_2$ and
$\mu_2 \mid \mu_0$, we obtain that
$$
\card(\supp\mu_0 \cap \supp\mu_1)=\card(\supp\mu_1 \cap \supp\mu_2)
=\card(\supp\mu_2 \cap \supp\mu_0)=1\,.
$$
Let $N=\card(\supp\mu_{0} \cap \supp\mu_{1} \cap\supp\mu_{2})$.
Clearly $N\le 1$; let us show that $N\ne 0$.
Suppose to the contrary that $N=0$.
Relations \eqref{eqvs in proper reminder case} apply to
$\mu_0$ and $\mu_1$, and without loss of generality we can assume that
$i' = j'$. For $\mu_1$ and $\mu_2$ the same relations
tell us that either $j' = k'$ or $j'' = k''$; since $N=0$, the first
case is impossible and we are left with $j'' = k''$.
Finally, we apply \eqref{eqvs in proper reminder case}
to $\mu_0$ and $\mu_2$ and see that either $i' = k'$ or $i'' = k''$,
which leads to a contradiction and completes the proof.
\end{proof}
\vskip5pt

\noindent The ${\mathbf S}(n)$
action in $\mathcal{E}^n(\mathbb{T}^2)$ induces an
${\mathbf S}(n)$ action in $L(\mathcal{E}^n(\mathbb{T}^2))$
defined by $(\sigma\lambda)(q)=\lambda(\sigma^{-1}q)$. 
The latter action, in turn, induces a simplicial ${\mathbf S}(n)$ action
in the complex $L_\vartriangle(\mathcal{E}^n(\mathbb{T}^2))$
which preserves dimension of simplices; our nearest goal
is to describe the orbits of this action.

\begin{Definition}\label{Def: normal form of simplices}
We define the following {\em normal forms} of simplices
of dimension $s>0$: 
$\Delta^s_\mathfrak{m}=
\{e_{\mathfrak{m};1,2},...,e_{\mathfrak{m};1,s+2}\}$,
$\nabla^s_\mathfrak{m}
=\{e_{\mathfrak{m};2,1},...,e_{\mathfrak{m};s+2,1}\}$,
where $\mathfrak{m}\in\mathfrak{M}_+$;
these simplices are called {\em normal}.
\end{Definition}

\begin{Lemma}\label{Lm: Sn action on simplicial complex}
For $s>0$, there are exactly $\card\mathfrak{M}$ orbits
of the ${\mathbf S}(n)$ action on the set of all $s$-simplices.\footnote{The 
number of ${\mathbf S}(n)$ orbits on the set of all $0$-simplices
is $\card\mathfrak{M}/2$.}
Every orbit contains exactly one normal simplex.
\end{Lemma}

\begin{proof}
Since $e_{\mathfrak{m};a,b} \nmid e_{\mathfrak{m};b,c}$,
Lemma \ref{lm: special property of the L complex} implies that
any $s$-simplex  $\Delta\in
L_{\vartriangle}({\mathcal C}^n(\mathbb{T}^2))$ is either of the form
$\{e_{\mathfrak{m};a,b_0},...,e_{\mathfrak{m};a,b_{s}}\}$
or of the form
$\{e_{\mathfrak{m};b_0,a},...,e_{\mathfrak{m};b_{s},a}\}$
with some $\mathfrak{m}\in\mathfrak{M}_+$ and distinct $a,b_0,...,b_s$.
An appropriate permutation $\sigma\in{\mathbf S}(n)$ carries
$\Delta$ to a normal form.
\end{proof}

\begin{Corollary}\label{Crl: dim of simplicial complex}
$\dim L_{\vartriangle}({\mathcal E}^n(\mathbb{T}^2))=n-2$.
\end{Corollary}

\subsection{Regular mappings
${\mathcal E}^n({\mathbb T}^2)\to{\mathbb T}^2$}
\label{ss: Regular mappings En(T2) to T2}

\noindent The following lemma will be used in Section
\ref{Sec: Holomorphic mappings of configuration spaces}.

\begin{Lemma}\label{lm: big dimensional torus mapping to torus}
Let
$\lambda:(\mathbb{T}^2)^n\to\mathbb{T}^2$
be a rational map. Then it is regular and
there are a holomorphic self-maps
$\mu_1,...,\mu_n$ of $\mathbb{T}^2$ such that
$\lambda(q_1,...,q_n)=\sum\limits_{i=1}^n\mu_i(q_i)$.
In particular,
$$
\lambda(q_1,...,q_n)=\sum\limits_{i=1}^n
\sum\limits_{\mathfrak{m}\in\mathfrak{M}_+}k_{i,\mathfrak{m}}
\mathfrak{m}q_i+c \ \ \text{\rm for all} \ 
(q_1,...,q_n)\in\mathcal{E}^n(\mathbb{T}^2)\,,
$$
where $k_{i,\mathfrak{m}} \in \mathbb{Z}$ and $c\in \mathbb{T}^2$.
\end{Lemma}
\begin{proof}
The proof is by induction on $n$. For $n=1$,
$\lambda$ is a rational map of smooth projective curves.
It extends to a regular map
(see, for instance, \cite{Shafarevich94}, Chapter II, Sec. 3.1, Corollary 1).
By Claim B, the regular self-map $\lambda$ of $\mathbb{T}^2$ is of
the desired form.
\vskip5pt

\noindent Assume that the theorem has already been proved
for some $n=m-1\ge 1$.
There is a subset $\Sigma\subset(\mathbb{T}^2)^m$ of codimension $1$
such that $\lambda$ is regular on $(\mathbb{T}^2)^m\setminus\Sigma$.
Let $(t_0,z_0)\in (\mathbb{T}^2\times(\mathbb{T}^2)^{m-1})\setminus \Sigma$
and $D$ be a small neighbourhood of $z_0$
in $(\mathbb{T}^2)^{m-1}$. Without loss of generality, we may assume that  
$t_0=0$ and $(0,z)\notin\Sigma$ for all $z\in D$. 
For $(t,z)\in (\mathbb{T}^2\times D)\setminus\Sigma$,
set $\mu(t,z)=\lambda(t,z)-\lambda(0,z)$ 
and $\nu(t,z)=\mu(t,z)-\mu(t,z_0)$. For any $z\in D$, we have
$\nu(0,z)=0$ and hence the map $t\mapsto\nu(t,z)$ extends to
a holomorphic endomorphism $\nu_z$ of the torus $\mathbb{T}^2$; 
moreover, the endomorphism $\nu_{z_0}$ carries the whole torus
$\mathbb{T}^2$ to the zero point $0\in\mathbb{T}^2$.
One can find a neighbourhood $D'\Subset D$ of $z_0$
and a compact subset $K\subset\mathbb{T}^2\times D$
such that for all $z\in D'$ the intersection
$K\cap (\mathbb{T}^2\times\{z\})$ is a union of two loops
that do not meet $\Sigma$ and generate the whole fundamental
group of the torus $\mathbb{T}^2\times\{z\}$. Moreover, since 
$\nu(\mathbb{T}^2\times\{z_0\})=0$, we may also assume that
$\nu(K)$ is contained in a small contractible neighbourhood of the zero point
in $\mathbb{T}^2$. It follows that for any $z\in D'$
the endomorphism $\nu_{z}$ is contractible and hence trivial.
Thus, for such $z$ we have $\mu(t,z)-\mu(t,z_0)\equiv 0$ and
$\lambda(t,z)\equiv\lambda(0,z)+\lambda(t,z_0)-\lambda(0,z_0)$.
By the uniqueness theorem, the latter identity holds true on the whole
of $(\mathbb{T}^2\times(\mathbb{T}^2)^{m-1})\setminus \Sigma$; 
the inductive hypothesis applies to $\lambda(0,z)$ and $\lambda(t,z_0)$,
which leads to the desired representation of $\lambda(t,z)$ and completes the
proof. 
\end{proof}


\section{Holomorphic mappings of configuration spaces}
\label{Sec: Holomorphic mappings of configuration spaces}

\noindent The main goal of this section is to prove the
classification theorem \ref{Thm: tame and orbit-like maps}
for holomorphic endomorphisms of torus configuration spaces.
We are proceeding according to the plan exposed in Introduction.

\subsection{Strictly equivariant mappings}
\label{ss: Strictly equivariant mappings}
A lifting $f$ of a continuous self-map $F$ of ${\mathcal C}^n(X)$
to the covering ${\mathcal E}^n(X)$
is equivariant, meaning that there is an endomorphism
$\alpha$ of the symmetric group ${\mathbf S}(n)$
such that $f(\sigma q)=\alpha(\sigma)f(q)$ for all
$q\in{\mathcal E}^n(X)$ and $\sigma\in{\mathbf S}(n)$.
We will see that a lifting of a {\em non-abelian} $F$
has a stronger property defined below.
\vskip5pt

\begin{Definition}\label{Def: strict equivariance}
A continuous map $f\colon{\mathcal E}^n(X)\to{\mathcal E}^n(X)$
is said to be {\em strictly equivariant} if there exists
an automorphism $\alpha$ of the group ${\mathbf S}(n)$
such that
\begin{equation}\label{eq: strict equivariance}
f(\sigma q)=\alpha(\sigma)f(q) \ \
\text{for all} \ q\in{\mathcal E}^n(X) \
\text{and} \ \sigma\in{\mathbf S}(n)\,.
\end{equation}
\end{Definition}

\begin{Theorem}\label{Thm: Lifting and Equivariance Theorem}
$a)$ For $n>4$ any non-abelian continuous map
$F\colon{\mathcal C}^n(\mathbb{T}^2)\to
{\mathcal C}^n(\mathbb{T}^2)$ admits a
continuous lifting
$f\colon{\mathcal E}^n(\mathbb{T}^2)\to
{\mathcal E}^n(\mathbb{T}^2)$ which fits in the diagram
\eqref{CD: lifting diagram}.

$b)$ For $n>4$ any continuous lifting
$f\colon{\mathcal E}^n(\mathbb{T}^2)\to
{\mathcal E}^n(\mathbb{T}^2)$
of a non-abelian continuous map
$F\colon{\mathcal C}^n(\mathbb{T}^2)\to
{\mathcal C}^n(\mathbb{T}^2)$
is strictly equivariant.
\end{Theorem}

\begin{proof}
In view of the covering mapping theorem, $(a)$ 
follows from Theorem \ref{Thm: algebraic lifting for torus}.
Let us prove $(b)$. The diagram \eqref{CD: lifting diagram} for
$f$ and $F$ leads to the algebraic commutative diagram
\begin{equation}\label{CD: algebraic lifting diagram}
\CD
1@>>>\pushvertex{-15pt}{\pi_1({\mathcal E}^n(\mathbb{T}^2),Q^\circ)}
          @>{p_*}>> \pushvertex{-15pt}{\pi_1({\mathcal C}^n(\mathbb{T}^2),Q^\circ)}
          @>{\delta}>>\pushvertex{-15pt}{\mathbf S}(n)\arrowscale{0.5}@>>> \pushvertex{-18pt}1\\
@. \nudge{-15pt} @V{f_*}VV \nudge{-15pt}@V{F_*}VV 
\nudge{-15pt}@ V{\alpha}VV \\ 
1@>>>\pushvertex{-15pt}{\pi_1({\mathcal E}^n(\mathbb{T}^2),f(Q^\circ))}
@>>{p_*}> \pushvertex{-15pt}{
\pi_1({\mathcal C}^n(\mathbb{T}^2),f(Q^\circ))} @>>{\delta}> 
\pushvertex{-15pt}{\mathbf S}(n)\arrowscale{0.5}@>>>
\pushvertex{-16pt} 1\,,
\endCD
\end{equation}
which relates the final segments of the exact homotopy sequences 
of the coverings $p$ in the left and right columns of
\eqref{CD: lifting diagram}.
The condition \eqref{eq: strict equivariance} holds true with 
the endomorphism $\alpha$ that appears in \eqref{CD: algebraic lifting diagram},
and we have only to show that this $\alpha$ is an automorphism whenever
$F^*$ is non-abelian and $n>4$. 

Suppose to the contrary that $\alpha$ is not an automorphism; then
its image is a non-trivial quotient of ${\mathbf S}(n)$, which
must be abelian since $n>4$. In view of \eqref{CD: algebraic lifting diagram},
the composition $\delta\circ F^*=\alpha\circ\delta$ is abelian and,
by Lemma \ref{Lm: composition of homomorphisms},
$F^*$ itself is abelian, which contradicts our assumption.  
\end{proof}

\noindent The following lemma shows that the assertion of Lemma
\ref{Lm: simplicial map f*: L(Y) to L(Z)} holds true
for every strictly equivariant map.

\begin{Lemma}\label{Lm: nonconstant differences}
Let $n>2$ and $f=(f_1,...,f_n)\colon{\mathcal E}^n(\mathbb{T}^2)
\to{\mathcal E}^n(\mathbb{T}^2)$ be a strictly equivariant
holomorphic map. Then
$$
f^*\colon L({\mathcal E}^n(\mathbb{T}^2))
\ni\lambda\mapsto\lambda\circ f\in L({\mathcal E}^n(\mathbb{T}^2))
$$
is a well-defined simplicial map whose restriction to each
simplex $\Delta\subseteq L({\mathcal E}^n(\mathbb{T}^2))$
is injective. Consequently,
$f^*$ preserves dimension of simplices.
\end{Lemma}

\begin{proof}
In view of Lemma \ref{Lm: simplicial map f*: L(Y) to L(Z)}, we
have only to prove that $\mu\circ f\ne\const$ for any
$\mu\in L({\mathcal E}^n(\mathbb{T}^2))$.
Suppose to the contrary that $\mu\circ f=c\in \mathbb{T}^2$.
Then $(\mu\circ f)(\sigma q)\equiv c$
for all $\sigma \in{\mathbf S}(n)$. Since $f$ is strictly
equivariant, there is $\alpha\in\Aut{\mathbf S}(n)$ such that
$f(\sigma q)=\alpha(\sigma) f(q)$ for all
$\sigma\in{\mathbf S}(n)$ and $q\in{\mathcal E}^n(X)$,
so that $c\equiv \mu(f(\sigma q))=\mu(\alpha(\sigma) f(q))$.

By Theorem \ref{Thm: description of all differences},
$\mu=\mathfrak{m}(q_i-q_j)$ for some
distinct $i,j$ and $\mathfrak{m}\in \mathfrak{M}$;
hence $c\equiv(\mu\circ f)(q)=\mathfrak{m}(f_i(q)-f_j(q))$.
Since $\alpha$ is an automorphism and $n>2$, there is
$\sigma\in{\mathbf S}(n)$ such that
$\alpha(\sigma^{-1})(i)=i$ and $\alpha(\sigma^{-1})(j)=k\ne j$;
thus, for all $q$ we have
$$
\aligned
\mathfrak{m}(f_i(q)-f_j(q))&=c=(\mu\circ f)(\sigma q)
=\mu(\alpha(\sigma) f(q))\\
&=\mathfrak{m}(f_{\alpha(\sigma^{-1})(i)}(q)
-f_{\alpha(\sigma^{-1})(j)}(q))=\mathfrak{m}(f_i(q)-f_k(q))\,,
\endaligned
$$
which is impossible.
\end{proof}

\subsection{Proof of Theorem \ref{Thm: tame and orbit-like maps}}
\label{ss: Proof of main theorem}
\noindent We shall prove two theorems, which 
together imply
Theorem \ref{Thm: tame and orbit-like maps}.

\begin{Theorem}\label{Thm: tame torus}
For $n>4$, every holomorphic non-abelian self-map $F$ of
$\mathcal{C}^n(\mathbb{T}^2)$ is tame.
\end{Theorem}
\begin{proof}
By Theorems \ref{Thm: algebraic lifting for torus} and
\ref{Thm: Lifting and Equivariance Theorem},
the map $F$ lifts to a strictly equivariant holomorphic
map $f$ that fits to the commutative diagram
\begin{equation*}
\CD
\mathcal{E}^n(\mathbb{T}^2)@>{f}>> \mathcal{E}^n(\mathbb{T}^2) \\
@V{p}VV @VV{p}V\\
\mathcal{C}^n(\mathbb{T}^2) @>{F}>>\mathcal{C}^n(\mathbb{T}^2)\,.
\endCD
\end{equation*}
Let $\alpha$ be the automorphism of ${\mathbf S}(n)$ corresponding
to $f$ (see Definition \ref{Def: strict equivariance}).

By Lemma \ref{Lm: nonconstant differences},
$f^*$ is a dimension preserving simplicial self-map of
$L_\Delta(\mathcal{E}^n(\mathbb{T}^2))$. Let
$\Delta_1 =\{q_1 - q_2,...,q_1 - q_n\}$ and
$\Delta=f^*(\Delta_1)$. By Lemma
\ref{Lm: Sn action on simplicial complex},
there exists a permutation
$\sigma$ that brings $\Delta$ to its normal form;
without loss of generality, we may assume that this normal form
is $\nabla_{\mathfrak m}
=\{\mathfrak{m}(q_2 - q_1),...,\mathfrak{m}(q_n - q_1)\}$,
where $\mathfrak{m}\in\mathfrak{M}_+$.
Set $\tilde{f}=f\circ \sigma$; then
\begin{equation}\label{eq: difference}
\tilde{f}_j=\tilde{f}_1 + \mathfrak{m}(q_1 - q_j)\,,  \ \
j=1,...,n\,.
\end{equation}
Define the holomorphic map
$\tau\colon \mathcal{E}^n(\mathbb{T}^2) \to \Aut(\mathbb{T}^2)$
by the condition
$$\tau(q)(z)=\tau(q_1,...,q_n)(z)
=-\mathfrak{m}z + (\tilde{f}_1(q) + \mathfrak{m}q_1)\,,
$$
where
$q=(q_1,...,q_n)\in \mathcal{E}^n(\mathbb{T}^2)$ and
$z\in\mathbb{T}^2$.
Equations
\eqref{eq: difference} imply that $\tau(q)q_j=f_j(\sigma q)$ for
all $j=1,...,n$ and $q=(q_1,...,q_n)
\in{\mathcal E}^n(\mathbb{T}^2)$;
thereby $\tau(q)q=f(\sigma q)=\alpha(\sigma)f(q)$,
or, which is the same,
$f(q)=\alpha(\sigma^{-1})\tau(q)q$ for all $q\in{\mathcal
E}^n(\mathbb{T}^2)$. To complete the proof, we must check that
$\tau$ is $\mathbf{S}(n)$-invariant.
For every $s\in \mathbf{ S}(n)$ and
all $q\in{\mathcal E}^n(\mathbb{T}^2)$ we have
\begin{equation*}
\aligned
s\tau(sq)q&=\tau(sq)sq=f(\sigma sq)=\alpha(\sigma s)f(q)\\
&=\alpha(\sigma s)f(\sigma\sigma^{-1}q)
=\alpha(\sigma s)\tau(\sigma^{-1}q)\sigma^{-1}q
=\alpha(\sigma s)\sigma^{-1}\tau(\sigma^{-1}q)q\,,
\endaligned
\end{equation*}
which can be rewritten as
\begin{equation}\label{eq: a}
[(\tau(sq))^{-1}\cdot\tau(\sigma^{-1}q)]q
=\sigma\alpha(s^{-1}\sigma^{-1})sq\,,
\end{equation}
where
$(\tau(sq))^{-1}\cdot\tau(\sigma^{-1}q)\in\Aut(\mathbb{T}^2)$ is
the product in group $\Aut(\mathbb{T}^2)$. Let us notice that for
$n>1$ there is a non-empty Zariski open subset
$A\subset{\mathcal E}^n(\mathbb{T}^2)$
such that if $\theta q=\rho q$ for some $q\in A$,
$\theta\in\Aut(\mathbb{T}^2)$ and $\rho\in \mathbf{S}(n)$, then
$\theta=\id$ and $\rho=1$. Therefore, equation (\ref{eq: a})
implies $\tau(sq)=\tau(\sigma^{-1}q)$ and
$\sigma\alpha(s^{-1}\sigma^{-1})s=1$ for all $q\in A$ and all
$s\in \mathbf{S}(n)$. Since $\tau$ is continuous, the first of
these relations holds true for all
$q\in {\mathcal E}^n(\mathbb{T}^2)$
and all $s\in \mathbf{S}(n)$, which shows that
$\tau\colon {\mathcal E}^n(\mathbb{T}^2)\to\Aut(\mathbb{T}^2)$ is
invariant\footnote{We have also proved that $\alpha(s)=\sigma^{-1}s\sigma$}.
\end{proof}
\begin{Remark}
\label{Rm: automorphisms of Cn(T2) for n=3 or 4 are tame}
Let $n$ be either $3$ or $4$. The statement of Theorem
\ref{Thm: tame torus} still holds true if
we assume that the map $F$ is an automorphism.
The only changes we need to make in the proof are as follows:
instead of our Theorem \ref{Thm: Lifting and Equivariance Theorem},
we have to use Theorem $2$ from \cite{Ivanov92},
which states that the pure braid group is a characteristic
subgroup of the torus braid group; moreover, instead of
Lemma \ref{Lm: nonconstant differences}, we should use Remark
\ref{Rmk: dominant map} from Section 
\ref{ss: Simplicial structure in L(En(T2))}. The rest of the
proof is the same.
\end{Remark}

\begin{Remark}\label{Rm: uniqueness of tame form}
$a)$ {\sl Let $n\ge 2$ and let 
$F\colon{\mathcal C}^n(\mathbb{T}^2)
\to{\mathcal C}^n(\mathbb{T}^2)$ be a tame map.
Then a morphism $T\colon{\mathcal C}^n(\mathbb{T}^2)
\to\Aut{\mathbb T}^2$ in the 'tame representation'
$F=F_T$ of $F$ is uniquely determined by $F$.}
Indeed, if $F_T=F_{T'}$ for two morphisms $T,T'$, then
$T(Q)Q=T'(Q)Q$ and $(*)$ $[T(Q)]^{-1}T'(Q)Q=Q$
for all $Q\in{\mathcal C}^n(\mathbb{T}^2)$. Furthermore,
a torus automorphism is uniquely determined by its values
at a generic pair of distinct points;
since $n\ge 2$, the identity $(*)$ shows that
$[T(Q)]^{-1}T'(Q)=\id$ for any generic point  
$Q\in{\mathcal C}^n(\mathbb{T}^2)$ and hence $T(Q)=T'(Q)$
everywhere.
\vskip5pt

\noindent $b)$ In view of Theorem \ref{Thm: tame torus},
$(a)$ shows that {\sl for $n>4$ any holomorphic non-abelian 
map $F\colon{\mathcal C}^n(\mathbb{T}^2)
\to{\mathcal C}^n(\mathbb{T}^2)$ admits a unique 
tame representation $F=F_T$ and the morphism
$T$ is regular whenever $F$ is so}. 
Remark \ref{Rm: automorphisms of Cn(T2) for n=3 or 4 are tame}
shows that the uniqueness (and regularity)
of $T$ still hold true whenever $n$ is $3$ or $4$
and $F$ is a (biregular) automorphism.  
\end{Remark}

\begin{Definition}\label{ref: def of m(z)}
The map
\begin{equation*}
{\mathcal C}^n(\mathbb{T}^2)\ni q=\{q_1,\dots,q_n\}
\mapsto {\mathit s}(q)= (q_1 +\dots+q_n) \in \mathbb{T}^2
\end{equation*}
is a locally trivial holomorphic fibring whose fibre
$M_0 = s^{-1}(0)$ is an algebraic variety.
We refer to the variety $M_0$ as the
{\em reduced configuration space}.
The presentation of the fundamental group 
$\pi_1(M_0)$ found by O. Zariski \cite{Zariski37}
shows that $H_1(M_0,\mathbb{Z}) = \mathbb{Z}_{2n}$;
Zariski called the group $\pi_1(M_0)$
the {\em invariant subgroup of the
group of motion classes of an elliptic Riemann surface}.
\end{Definition}

\noindent Let $\gamma:\mathbb{C} \to \mathbb{T}^2$
be the universal covering; then the mappings
\begin{equation}\label{eq: coverings maps}
\aligned
&M_0\times \mathbb{C} \ni (q,\zeta) \mapsto {{\mathit h}}(q,\zeta)
=\{q_1+\gamma(\zeta),..., q_n+\gamma(\zeta)\}
\in {\mathcal C}^n(\mathbb{T}^2)\,,\\
&M_0\times \mathbb{T}^2 \ni (q,t) \mapsto \tilde{{\mathit h}}(q,t)
= \{q_1+t,..., q_n+t\} \in {\mathcal C}^n(\mathbb{T}^2)
\endaligned
\end{equation}
are holomorphic coverings, as well.
\vskip5pt

\noindent The following theorem completes the classification of
endomorphisms of torus configuration spaces.

\begin{Theorem}\label{orbit-like}
If $m>2$, then a holomorphic map
$F\colon\mathcal{C}^n(\mathbb{T}^2)\to
\mathcal{C}^m(\mathbb{T}^2)$ is orbit-like if and only if
it is abelian.
\end{Theorem}

\begin{proof}
Let $F$ be abelian. By the abelianization of defining
relations \eqref{rel:Zariski 1} - \eqref{rel:Zariski 6},
we see that
$H_1({\mathcal C}^n(\mathbb{T}^2),\mathbb{Z})
= B_n(\mathbb{T}^2)/B_n^\prime(\mathbb{T}^2)=
\mathbb{Z}_2 \oplus \mathbb{Z}^2$.
As it was already mentioned,
$B_m(\mathbb{T}^2)=\pi_1({\mathcal C}^m(\mathbb{T}^2))$
has no torsion. It follows that
$\Img F_*$ is either $\mathbb{Z}^2$ or $\mathbb{Z}$ or trivial.
Since $\pi_1(M_0)/(\pi_1(M_0))^\prime = \mathbb{Z}_{2n}$, there is
no non-trivial homomorphism $\pi_1(M_0) \to \Img F_*$; that is,
the homomorphism $(F\circ h)_*$ is trivial.
This implies that there exists a holomorphic map $f=(f_1,...,f_m)$
which fits to the commutative diagram
\begin{equation*}
\CD
M_0\times \mathbb{C}@>{f}>>{\mathcal E}^m(\mathbb{T}^2)\\
@V{\mathit h}VV@/SE/ F\circ h//@VV{p}V\\
{\mathcal C}^n(\mathbb{T}^2) @>>{F}>
{\mathcal C}^m(\mathbb{T}^2)\,,
\endCD
\end{equation*}
where $h$ is defined in	\eqref{eq: coverings maps}.
The homomorphism $f_*$ of the fundamental groups
induced by $f$ is trivial. Hence, for any $j$,
the composition
$$
f_j-f_1=(q_j-q_1)\circ f
\colon M_0\times \mathbb{C}\stackrel{f}
{\longrightarrow}{\mathcal E}^m(\mathbb{T}^2)
\stackrel{q_j-q_1}{\longrightarrow}\mathbb{T}^2\setminus\{0\}
$$ 
is contractible and lifts to
a holomorphic map $g_j\colon M_0\times\mathbb{C}\to\mathbb{D}
=\{\zeta\in{\mathbb C}\,|\ |\zeta|<1\}$ into the universal covering $\mathbb{D}$
of $\mathbb{T}^2\setminus\{0\}$.
Since $M_0\times\mathbb{C}$ is algebraic, the Liouville's theorem
shows that $g_j=\const$ and, thereby, $f_j-f_1=\const=c_j
\in\mathbb{T}^2\setminus\{0\}$.  
Thus, $f(q) = (0+f_1(q),c_2+f_1(q),\dots,c_m+f_1(q))$,
which shows that $f$ is orbit-like.
\vskip5pt

\noindent Suppose now that $F$ is orbit-like. To prove that
$F$ is abelian, it suffices to show that
for any point $q\in{\mathcal C}^m(\mathbb{T}^2)$,
the fundamental group of any connected component
of the $(\Aut\mathbb{T}^2)$-orbit ${\mathcal O}_q
=(\Aut\mathbb{T}^2)(q)$ is abelian.
For $m>2$, any component of ${\mathcal O}_q$ is diffeomorphic to the orbit
${\mathcal O}_q^*$ of the action of $\mathbb{T}^2$ in
${\mathcal C}^m(\mathbb{T}^2)$ by translations. The latter orbit
${\mathcal O}_q^*$ is a quotient group of $\mathbb{T}^2$
by a finite subgroup and hence is homeomorphic to $\mathbb{T}^2$.
Thus, $\pi_1({\mathcal O}_q^*)={\mathbb Z}^2$.
\end{proof}

\noindent We conclude this section with two simple
statements about abelian maps. 

\subsection{Splitting of abelian maps}
\label{ss: Splitting of abelian maps}
Up to a homotopy, any abelian map
$$
f\colon\mathcal{C}^n(\mathbb{T}^2)\to
\mathcal{C}^m(\mathbb{T}^2)
$$
splits to a composition
$g\circ s$ of the standard map 
$s\colon\mathcal{C}^n(\mathbb{T}^2)\to\mathbb{T}^2$,
$s(q)=q_1+...+q_n$, defined in Definition
\ref{ref: def of m(z)} and an appropriate continuous map 
$g\colon\mathbb{T}^2\to\mathcal{C}^m(\mathbb{T}^2)$.
Thus, we obtain the commutative up to a homotopy
diagram
$$
\CD
{\mathcal C}^n(\mathbb{T}^2)@[2]>{f}>>{\mathcal C}^m(\mathbb{T}^2)\\
 @. @/SE//{{\mathit s}}/@. @/NE/{g}//\\
  @. \pushvertex{4pt}\mathbb{T}^2\,.
\endCD
$$
Indeed, let $p\colon\mathcal{E}^n(\mathbb{T}^2)\to
\mathcal{C}^n(\mathbb{T}^2)$ be the standard projection
and  
$$
\gamma=(\gamma(t),a_2,...,a_n)
\subset\mathcal{E}^n(\mathbb{T}^2)
$$
be a loop.
The map $s\circ p$ carries $\gamma$ to the loop
$\gamma(t)+a_2+...+a_n$ in $\mathbb{T}^2$, which shows that 
$(s\circ p)_*$ is an epimorphism. Therefore, $s_*$ is an
epimorphism as well.
The homomorphism $s_*$ splits to a composition
$s_*=\phi\circ\alpha$,
where $\alpha\colon B_n(\mathbb{T}^2)
\to B_n(\mathbb{T}^2)/[B_n(\mathbb{T}^2),B_n(\mathbb{T}^2)]
=\mathbb{Z}^2\times\mathbb{Z}_2\to\mathbb{Z}^2$
is the composition of the abelianization 
and the torsion eliminating map 
$\mathbb{Z}^2\times\mathbb{Z}_2\to\mathbb{Z}^2$,
and $\phi\colon \mathbb{Z}^2\to\mathbb{Z}^2$
is an epimorphism. Since every surjective endomorphism of
$\mathbb{Z}^2$ is an isomorphism,
it follows that $\phi^{-1}\circ s_* =\alpha$.
Since the map $f$ is abelian there exists a homomorphism
$\beta\colon\mathbb{Z}^2 \to B_m(\mathbb{T}^2)$ such that
$f_*=\beta \circ \alpha$.
$\mathcal{C}^m(\mathbb{T}^2)$ is a $K(\pi,1)$ space
for $\pi=B_m(\mathbb{T}^2)$, which
implies that there is a continuous map
$g\colon\mathbb{T}^2\to\mathcal{C}^m(\mathbb{T}^2)$ such that
$g_*=\beta \circ \phi^{-1}$. Clearly
$f$ is homotopic to $g\circ s$, which proves our statement.

\subsection{Holomorphic mappings $\mathbb{T}^2\to
\mathcal{C}^n(\mathbb{T}^2)$}
\label{ss: Holomorphic mappings T2 to Cn(T2)}

\begin{Proposition}\label{prop: T to configurations is orbit-like}
Any holomorphic map 
$F\colon\mathbb{T}^2\to \mathcal{C}^n(\mathbb{T}^2)$
carries $\mathbb{T}^2$ to an orbit of
the $\Aut \mathbb{T}^2$ action in
$\mathcal{C}^n(\mathbb{T}^2)$.
\end{Proposition}
\begin{proof}
$F$ lifts to a holomorphic map $f$ of the universal coverings
and we have the commutative diagram
$$
\CD
\mathbb{C}@>{f=(f_1,...,f_n)}>>
{\mathcal E}^n(\mathbb{T}^2)\\
@V{{\mathit \gamma}}VV @VV{p}V\\
\mathbb{T}^2 @>>{F}>{\mathcal C}^n(\mathbb{T}^2)\,.
\endCD
$$
Each $f_j$ maps $\mathbb{C}$ to $\mathbb{T}^2$ and,
for $j>1$, $f_j-f_1$ maps $\mathbb{C}$ to 
$\mathbb{T}^2\setminus\{0\}$. It follows that 
$f_j-f_1=\const=c_j\in\mathbb{T}^2\setminus\{0\}$
and $f = (0+f_1,c_2+f_1,\dots,c_n+f_1)$.
\end{proof}

\section{Biregular automorphisms}
\label{Sec: Biregular automorphisms}

\noindent Here we describe all biregular automorphisms of
the algebraic manifold $\mathcal{C}^n(\mathbb{T}^2)$.

\begin{Lemma}\label{Lm: c^n to T}
Any regular map
$R\colon\mathcal{C}^n(\mathbb{T}^2)\to \mathbb{T}^2$
is of the form
\begin{equation*}
R(Q)=\sum\limits_{\mathfrak{m}\in \mathfrak{M}_+}
k_{\mathfrak{m}}\mathfrak{m}(q_1+...+q_n) + c \ 
\text{for all} \ Q=\{q_1,...,q_n\} \in 
\mathcal{C}^n(\mathbb{T}^2)\,,
\end{equation*}
where $k_{\mathfrak{m}} \in \mathbb{Z}$ and $c\in \mathbb{T}^2$.
\end{Lemma}

\begin{proof}
Consider the map $r=R\circ p$, where
$p\colon \mathcal{E}^n(\mathbb{T}^2)\to
\mathcal{C}^n(\mathbb{T}^2)$ is the standard projection.
By Lemma \ref{lm: big dimensional torus mapping to torus},
$r(q)=\sum\limits_{i=1}^n
\sum\limits_{\mathfrak{m}\in \mathfrak{M}_+}k_{i,\mathfrak{m}}
\mathfrak{m}q_i + c$.
Since $r$ must be invariant under the $\mathbf{S}(n)$ action,
it follows that $k_{1,\mathfrak{m}}=...=k_{n,\mathfrak{m}}=k_{\mathfrak{m}}$.
Thus,
$$
r(q_1,\dots,q_n)=\sum\limits_{i=1}^n
\sum\limits_{\mathfrak{m} \in
\mathfrak{M}_+}k_{\mathfrak{m}}\mathfrak{m}q_i +
c=\sum\limits_{\mathfrak{m}\in
\mathfrak{M}_+}k_{\mathfrak{m}}\mathfrak{m}(q_1+...+q_n) + c
$$
and $R(Q)=\sum\limits_{\mathfrak{m}\in
\mathfrak{M}_+}k_{\mathfrak{m}}\mathfrak{m}(q_1+...+q_n) + c$.
\end{proof}

\begin{Theorem}\label{Aut}
For $n>2$, any biregular automorphism $F$ of
${\mathcal C}^n(\mathbb{T}^2)$ is of the form $F(Q)=AQ$,
where $A \in \Aut\mathbb{T}^2$.
\end{Theorem}

\begin{proof}
\noindent Since for $n>2$ the group $B_n(\mathbb{T}^2)$
is non-abelian, $F$ is non-abelian.
By Theorem \ref{Thm: tame and orbit-like maps} and
Remark \ref{Rm: uniqueness of tame form} from Section
\ref{ss: Proof of main theorem},
there is a unique regular map
$T\colon {\mathcal C}^n(\mathbb{T}^2)\to \Aut \mathbb{T}^2$ 
such that $F(Q)=T(Q)Q$ for all
$Q\in {\mathcal C}^n(\mathbb{T}^2)$.

For any $Q\in{\mathcal C}^n(\mathbb{T}^2)$,
$T(Q)$ is an automorphism of
$\mathbb{T}^2$ and, by Claim B in Section
\ref{ss: Holomorphic mappings En to T2-0} and 
Lemma \ref{Lm: c^n to T}, 
it carries a point $z\in\mathbb{T}^2$ to the point
\begin{equation}\label{eq: 5.1}
T(Q)z= \mathfrak{m}_0 z +
\sum\limits_{\mathfrak{m}\in \mathfrak{M}_+}
k_{\mathfrak{m}}\mathfrak{m}(q_1+...+q_n) + c\,, 
\end{equation}
where $\mathfrak{m}_0\in\mathfrak{M}$, 
$k_{\mathfrak{m}} \in \mathbb{Z}$, and $c\in \mathbb{T}^2$
do not depend on $z$ and $Q$.

Recall that $s(Q)=q_1+...+q_n$
for $Q=\{q_1,...,q_n\}$ and set $s_1=s\circ F$, that is,
$$
s_1(Q)=(s\circ F) (Q) = s(T(Q)Q) =
T(Q)q_1+...+T(Q)q_n\,.
$$
Using \eqref{eq: 5.1} for $z=q_1,...,q_n$, we see that
\begin{equation}\label{eq: 5.2}
\aligned
{\mathit s}_1(Q)&=\mathfrak{m}_0(q_1+...+q_n)+
n\big(\sum\limits_{\mathfrak{m}\in\mathfrak{M}_+}
             k_{\mathfrak{m}}\mathfrak{m}(q_1+...+q_n) + c\big)\\
&=\big(\mathfrak{m}_0
+n\sum\limits_{\mathfrak{m}\in\mathfrak{M}_+}k_{\mathfrak{m}}
\mathfrak{m}\big)(q_1+...+q_n) + n c\,.
\endaligned
\end{equation}

On the other hand, $F^{-1}$ is a regular automorphism as well.
By the same argument, there is a unique regular
$T^{\prime}\colon{\mathcal C}^n(\mathbb{T}^2)\to
\Aut \mathbb{T}^2$ such that $F^{-1}(Q)=T^{\prime}(Q)Q$ for
$Q\in {\mathcal C}^n(\mathbb{T}^2)$. As above, we conclude
that $T^{\prime}(Q)$ carries a point $z\in\mathbb{T}^2$
to the point
\begin{equation}\label{eq: 5.3}
T'(Q)z= \mathfrak{m}_0' z +
\sum\limits_{\mathfrak{m}\in \mathfrak{M}_+}
k_{\mathfrak{m}}'\mathfrak{m}(q_1+...+q_n) + c'\,, 
\end{equation}
where $\mathfrak{m}_0'\in\mathfrak{M}$, 
$k_{\mathfrak{m}}' \in \mathbb{Z}$, and $c'\in \mathbb{T}^2$
do not depend on $z$ and $Q$.
Since $s_1 \circ F^{-1}=s$, formulas
\eqref{eq: 5.2}--\eqref{eq: 5.3} imply
\begin{equation*}
\aligned
q_1+&...+q_n={\mathit s}(Q)= {\mathit s}_1(F^{-1}(Q))= 
{\mathit s}_1(\{T'(Q)q_1,...,T'(Q)q_n\})\\
&= \big(\mathfrak{m}_0
+n\sum\limits_{\mathfrak{m}\in\mathfrak{M}_+}k_{\mathfrak{m}}
\mathfrak{m}\big)\big(T'(Q)q_1+...+T'(Q)q_n\big) + n c \\
&=
\big(\mathfrak{m}_0
+n\sum\limits_{\mathfrak{m}\in\mathfrak{M}_+}k_{\mathfrak{m}}
\mathfrak{m}\big)\sum\limits_{i=1}^n\big(\mathfrak{m}_0' q_i +
\sum\limits_{\mathfrak{m}\in \mathfrak{M}_+}
k_{\mathfrak{m}}'\mathfrak{m}(q_1+...+q_n) + c'\big) + n c\,.
\endaligned
\end{equation*}
By changing the order of the summation, we obtain

\begin{equation}\label{eq: composition}
\aligned
q_1+...+q_n  =& \big(
\mathfrak{m}_0+n\sum\limits_{\mathfrak{m}\in \mathfrak{M}_+}
k_{\mathfrak{m}}\mathfrak{m} \big) \big(\mathfrak{m}_0'+ n\sum\limits_{\mathfrak{m}\in \mathfrak{M}_+}
k_{\mathfrak{m}}'\mathfrak{m}\big) (q_1+...+q_n) \\
&+  n(\mathfrak{m}_0+n\sum\limits_{\mathfrak{m}}
k_{\mathfrak{m}}\mathfrak{m}
 )c^\prime +nc\,.
 \endaligned
\end{equation}
Since $q_1+...+q_n$ runs over the whole torus, 
\eqref{eq: composition} shows that
the composition $\lambda=\mu\circ\nu=\nu\circ\mu$
of the torus endomorphisms
$\mu\colon z\mapsto \big(\mathfrak{m}_0 
+n\sum\limits_{\mathfrak{m}\in \mathfrak{M}_+}
k_{\mathfrak{m}}\mathfrak{m}\big) z$ and
$\nu\colon z\mapsto\big(\mathfrak{m}_0' 
+ n\sum\limits_{\mathfrak{m}\in \mathfrak{M}_+}
k_{\mathfrak{m}}^\prime\mathfrak{m}\big) z$
is the identity automorphism. Hence $\mu$ and $\nu$
are torus automorphisms sending $0$ to $0$ and,
by Claim B, $\mu(z)\equiv\mathfrak{m}_1 z$ with some 
$\mathfrak{m}_1 \in \mathfrak{M}$.
Consequently,
$\big(\mathfrak{m}_0 - \mathfrak{m}_1
+n\sum\limits_{\mathfrak{m}\in \mathfrak{M}_+}
k_{\mathfrak{m}}\mathfrak{m}\big) z\equiv 0$, i.e.
$\mathfrak{m}_0 - \mathfrak{m}_1
+n\sum\limits_{\mathfrak{m}\in \mathfrak{M}_+}
k_{\mathfrak{m}}\mathfrak{m}=0$. Since $n>2$ and elements of 
$\mathfrak{M}_+$ are linearly independent over $\mathbb Q$,
the latter cannot happen unless
$\sum\limits_{\mathfrak{m}\in
\mathfrak{M}_+}k_{\mathfrak{m}}\mathfrak{m}=0$.
Therefore, $T(Q)z=\mathfrak{m}_0 z+c$ for all
$Q$ and $z$, which completes the proof.
\end{proof}

\section{Configuration spaces of universal families}
Here we construct configuration spaces of
the universal Teichm{\" u}ller family of tori and describe
their holomorphic self-maps.

The Teichm{\" u}ller space $T(1,1)$
of tori with one marked point
is isomorphic to the upper half plane ${\mathbb H}^+$.
The group $H=\mathbb{Z}\times\mathbb{Z}$ acts discontinuously
and freely in the space $\mathcal{V}=
T(1,1)\times\mathbb{C}=\mathbb{H}^{+}\times\mathbb{C}$
by weighted translations
$(\tau,z)\mapsto (\tau,z+l+m\tau)$, $(l,m)\in H$.
Let $V(1,1)=\mathcal{V}/H$;
the map $\psi\colon\mathcal{V}\to V(1,1)$ is a covering
and the holomorphic projection
$\pi\colon V(1,1)\to\mathbb{H}^{+}=T(1,1)$ 
is called the {\em universal Teichm{\" u}ller family}
of tori with one marked point (see \cite{EarleKra74}, Sec. 4.11).
All fibres $\pi^{-1}(\tau)$ are tori; each of them carries
a natural group structure, marked points are
neutral elements and form a holomorphic section of $\pi$.

\begin{Definition}\label{Def: fibred configuration space}
Let $\mathcal{C}^n_\pi(V(1,1))$ be the complex subspace of the
configuration space $\mathcal{C}^n(V(1,1))$ of $V(1,1)$ consisting
of all $Q=\{q_1,...,q_n\}\in \mathcal{C}^n(V(1,1))$ such that
$\pi(q_1)=...=\pi(q_n)$. Define the holomorphic
projection $\rho\colon \mathcal{C}^n_\pi(V(1,1))\to T(1,1)$
by $\rho(Q)=\pi(q_1)=...=\pi(q_n)$,
$Q=\{q_1,...,q_n\}\in \mathcal{C}^n_\pi(V(1,1))$; 
the triple $\{\rho,\mathcal{C}^n_\pi(V(1,1)),T(1,1)\}$,
or simply $\rho\colon \mathcal{C}^n_\pi(V(1,1))\to T(1,1)$,
is said to be the {\em fibred configuration space
of the universal Teichm{\" u}ller family 
$\pi\colon V(1,1)\to T(1,1)$}.
A {\em fibred morphism} of fibred configuration spaces is 
a holomorphic map
$F\colon \mathcal{C}^n_\pi(V(1,1))
\to \mathcal{C}^m_\pi(V(1,1))$
which fits into the commutative diagram
\begin{equation}\label{special map condition}
\CD
\mathcal{C}^n_\pi(V(1,1)) @[2]>{F}>> \mathcal{C}^m_\pi(V(1,1))\\
@. @/SE//{\rho}/ @. @/SW/{\rho}// \\
@. \pushvertex{4pt}T(1,1)\,.
\endCD
\end{equation}
One may similarly define the corresponding {\em ordered}
fibred configuration spaces
${\theta}\colon \mathcal{E}^n_\pi(V(1,1))\to T(1,1)$
and their fibred morphisms.
The {\em fibred power} 
$$
{\theta}\colon (V(1,1))^n_\pi\to T(1,1)
$$
is defined by $(V(1,1))^n_\pi
=\{(q_1,...,q_n)\in(V(1,1))^n\,|\ \pi(q_1)=...=\pi(q_n)\}$
and ${\theta}(q_1,...,q_n)=\pi(q_1)$.
The natural covering map $\mu\colon\mathcal{E}^n_\pi(V(1,1))\to
\mathcal{C}^n_\pi(V(1,1))$ is also a fibred morphism.
The 
closed complex subspace
$S=\{(q_1,...,q_n)\in (V(1,1))^n_\pi \,| \ 
q_1=...=q_n\}$ of $(V(1,1))^n_\pi$ is called the {\em small
diagonal} of the fibred power
${\theta}\colon (V(1,1))^n_\pi\to T(1,1)$; it is
isomorphic to $V(1,1)$.
\end{Definition}

\begin{Remark}
\label{Rmk: fibred configuration spaces are irreducible}
The map 
$$
{{\Psi}}=(\psi,...,\psi)
\colon(\mathbb{H}^{+}\times\mathbb{C})^n\to (V(1,1))^n
$$
is a universal holomorphic covering and the preimage
$\mathcal{P}_n={{\Psi}}^{-1}((V(1,1))^n_\pi)
=\{((\tau_1,z_1),...,(\tau_n,z_n))
\in(\mathbb{H}^{+}\times\mathbb{C})^n
\,|\ \tau_1=...=\tau_n \}$ of $(V(1,1))^n_\pi\subset(V(1,1))^n$
is isomorphic to $\mathbb{H}^+\times \mathbb{C}^n$, i.e.,
it is irreducible and non-singular.
Since ${{\Psi}}$ is locally biholomorphic, the space
$(V(1,1))^n_\pi={{\Psi}}(\mathcal{P}_n)$ is
irreducible and non-singular and the map
${{\Psi}}|_{\mathcal{P}_n}\colon\mathcal{P}_n\to
(V(1,1))^n_\pi$ is a universal covering.
The space $\mathcal{E}^n_\pi(V(1,1))\subset(V(1,1))^n_\pi$
is the complement of a proper analytic subset of
$(V(1,1))^n_\pi$; hence $\mathcal{E}^n_\pi(V(1,1))$
is a connected complex manifold. Since 
$\mu\colon\mathcal{E}^n_\pi(V(1,1))\to\mathcal{C}^n_\pi(V(1,1))$ 
is a holomorphic covering, $\mathcal{C}^n_\pi(V(1,1))$
also is a connected complex manifold.

Notice that the preimage ${{{\mathcal F}}}=
{{\Psi}}^{-1}(\mathcal{E}^n_\pi(V(1,1)))$
of $\mathcal{E}^n_\pi(V(1,1))\subset(V(1,1))^n_\pi$
is an open dense subset of
$\mathcal{P}_n\cong\mathbb{H}^+\times \mathbb{C}^n$
and the maps ${{\Psi}}|_{{{\mathcal F}}}
\colon{{\mathcal F}}\to \mathcal{E}^n_\pi(V(1,1))$ and
$\mu\circ{{\Psi}}|_{{{\mathcal F}}}
\colon{{\mathcal F}}\to \mathcal{C}^n_\pi(V(1,1))$ are
holomorphic coverings.
\end{Remark}

\begin{Definition}\label{Def: special fibred morphisms}
Let $g\colon\mathcal{C}^n_\pi(V(1,1))\to V(1,1)$ be a
fibred morphism. Any point $Q\in \mathcal{C}^n_\pi(V(1,1))$ 
belongs to a certain fibre $\rho^{-1}(\tau)$, which is 
the configuration space $\mathcal{C}^n(\pi^{-1}(\tau))$
of the torus ${\mathbb T}^2_\tau=\pi^{-1}(\tau)$;
so $Q$ may be viewed as an $n$-point subset of ${\mathbb T}^2_\tau$.
Since $g$ is a fibred morphism, $g(Q)$ is a point of the same torus
${\mathbb T}^2_\tau$; thus, $Q+g(Q)$ and $-Q+g(Q)$ are
well-defined $n$-point subsets of ${\mathbb T}^2_\tau$, or,
which is the same, points of $\mathcal{C}^n({\mathbb T}^2_\tau)
\subset\mathcal{C}^n_\pi(V(1,1))$. This provides us with two
fibred maps $G_{\pm}
=\pm\Id+g\colon\mathcal{C}^n_\pi(V(1,1))\to
\mathcal{C}^n_\pi(V(1,1))$ defined by $Q\mapsto \pm Q+g(Q)$.
\end{Definition}

\begin{Lemma}
The fibred maps $G_{\pm}$ are holomorphic.
\end{Lemma}

\begin{proof}
There is a commutative diagram
\begin{equation}\label{eq: lifting of fibred map}
\CD
{{\mathcal U}}\\
@VV{\lambda}V @[2]/SE//{h}/\\ 
 {{\mathcal F}}  @[2].  T(1,1)\times\mathbb{C}\\
@VV{{\varphi}=\mu\circ{{\Psi}}}V @. @V{\psi}VV \\
\mathcal{C}^n_\pi(V(1,1)) @[2]>{g}>> V(1,1)\\
@. @/SE//{\rho}/ @. @/SW/{\pi}// \\
@. \pushvertex{4pt}T(1,1)\,,
\endCD
\end{equation}
where $\lambda\colon{\mathcal U}\to {\mathcal F}$
is the pull-back of the covering
$\psi\colon T(1,1)\times\mathbb{C} \to  V(1,1)$ along the map 
$g\circ {\varphi}\colon {\mathcal F} \to V(1,1)$
and $r$, as usual, is the restriction to
${\mathcal U}\subset 
{\mathcal F}\times T(1,1)\times\mathbb{C}$
of the projection
${\mathcal F}\times T(1,1)\times\mathbb{C}
\to T(1,1)\times\mathbb{C}$.
By the construction, $h$ is a fibred morphism.
Since ${\mathcal F}$ is a domain in $T(1,1)\times\mathbb{C}^n$,
it follows that $\lambda(u)=(\rho\circ{\varphi}\circ\lambda(u),
z_1(u),...,z_n(u))$ for $u\in {{\mathcal U}}$, where for all
$i=1,...,n$ the functions
$z_i\colon{{\mathcal U}}\to \mathbb{C}$ are holomorphic.
Denote the natural projection $T(1,1)\times\mathbb{C}\to
\mathbb{C}$ by $z$. Define the map 
$R={\Id+h}\colon{{\mathcal U}}\to\mathcal{P}_n$
by $R(q)=(\rho\circ{\varphi}\circ\lambda(u),
z_1(u)+z(h(u)),...,z_n(u)+z(h(u)))$ for 
$q\in {{\mathcal U}}$.
Since $\psi(\rho\circ{\varphi}\circ\lambda(u),
z_i(u))\neq \psi(\rho\circ{\varphi}\circ\lambda(u),
z_j(u))$
for any $u\in {{\mathcal U}}$ and $i\ne j$,
it follows that for all $m,n\in \mathbb{Z}$ we have
$z_i(u)-z_j(u)\neq m+n\rho(\varphi(\lambda(u)))$ and
$(z_i(u)+z(h(u)))-(z_j(u)+z(h(u)))
\neq m+n\rho(\varphi(\lambda(u)))$;
thus $R(u)\in {\mathcal F}$, that is,
$R({{\mathcal U}})\subset{{\mathcal F}}$.
Set  $R'={\varphi}\circ
R\colon{{\mathcal U}}\to\mathcal{C}^n_\pi(V(1,1))$.

Let us show that
for $u\in({\varphi}\circ\lambda)^{-1}(Q)$ we have
$R'(u)=Q+g(Q)$ and the map $R'$
induces the well-defined single-valued fibred morphism 
which is equal to
$G_{+}\colon \mathcal{C}^n_\pi(V(1,1))\to 
\mathcal{C}^n_\pi(V(1,1))$.
Indeed, for $u\in({\varphi}\circ\lambda)^{-1}(Q)$
we have 
$R(u)=(\rho(Q),z_1(u)+z(h(u)),...,z_n(u)+z(h(u)))$.
Notice that for $\tau\in T(1,1)$ and $x,y\in
\mathbb{C}$ we have $\psi(\tau,x+y)=\psi(\tau,x)+\psi(\tau,y)$,
where the latter summation is the sum in the torus
$\mathbb{T}^2_\tau=\psi^{-1}(\tau)$.
Let $\tau=\rho(Q)$ and $\zeta=z(h(u))$.
Therefore 
$$
\aligned
R'(u)&={\varphi}(R(u))=\mu({\Psi}(R(u)))
=\mu({\Psi}(\tau,z_1(u)+\zeta,...,z_n(u)+\zeta))\\
&=\mu(\psi(\tau,z_1(u)+\zeta),...,\psi(\tau,z_n(u)+\zeta))
\\
&=
\mu(
\psi(\tau,z_1(u))+\psi(\tau,\zeta),...,
\psi(\tau,z_n(u))+\psi(\tau,\zeta))
\\
&=
\mu(
\psi(\tau,z_1(u)),...,\psi(\tau,z_n(u)))+\psi(\tau,\zeta)
\\
&=
\mu({\Psi}((\tau,z_1(u)),...,(\tau,z_n(u)))+\psi(\tau,\zeta))\\
&=
\varphi(\tau,z_1(u),...,z_n(u))+\psi(\tau,\zeta))
= Q+g(Q)\,.
\endaligned
$$
The same argument shows that $G_{-}=-\Id+g$ is also holomorphic.
\end{proof}

\begin{Remark}
Similarly to the Definition \ref{Def: special fibred morphisms},
for any fibred morphism
$$
f=(f_1,...,f_n)\colon\mathcal{E}^n_\pi(V(1,1))\to
\mathcal{E}^n_\pi(V(1,1))
$$
one can 
define fibred morphisms 
$$
\Id+f\colon\mathcal{E}^n_\pi(V(1,1))
\to(V(1,1))^n_\pi
$$
and
$$
-\Id+f\colon\mathcal{E}^n_\pi(V(1,1))
\to(V(1,1))^n_\pi
$$
which take 
any point $q=(q_1,...,q_n)\in\mathcal{E}^n_\pi(V(1,1))$ 
respectively to 
$$
q+f(q)=(q_1+f_1(q),...,q_n+f_n(q))\in (V(1,1))^n_\pi
$$
and to 
$$
-q+f(q)=(-q_1+f_1(q),...,-q_n+f_n(q))\in (V(1,1))^n_\pi\,.
$$
\end{Remark}

\noindent The following theorem is an analogue of Theorem
\ref{Thm: tame and orbit-like maps} for the case of the
fibred morphisms.

\begin{Theorem}\label{Thm: tame tori 2^nd}
Let $n>4$. For any fibred non-abelian morphism
$F\colon\mathcal{C}^n_\pi(V(1,1))\to
\mathcal{C}^n_\pi(V(1,1))$, there is a fibred morphism
$g\colon\mathcal{C}^n_\pi(V(1,1))\to V(1,1)$
such that $F$ is either $\Id +g$ or
$-\Id +g$.
\end{Theorem}

\begin{proof}
According to Theorem \ref{Thm: tame and orbit-like maps},
for any $\tau\in T(1,1)$ there exists a unique holomorphic map
$T_\tau\colon\rho^{-1}(\tau)\to \Aut \pi^{-1}(\tau)$
such that $F(Q)=T_\tau(Q)Q$ for any
$Q\in\rho^{-1}(\tau)\subset \mathcal{C}^n_\pi(V(1,1))$. 
There are no complex multiplication on a generic torus.
Thus, for any generic $\tau\in T(1,1)$ and any
$Q\in\rho^{-1}(\tau)$, there exists $c_\tau(Q)\in\pi^{-1}(\tau)$
such that the automorphism
$T_\tau(Q)$ maps a point $z\in\pi^{-1}(\tau)$ either to
$z+c_\tau(Q)$ or to $-z+c_\tau(Q)$; 
notice that the representation is unique.
By Baire Category Theorem, there exists an open set
$D\subset T(1,1)$ and its dense subset $D'\subset D$
such that either $(*)$ for all $\tau \in D'$ we have
$T_\tau(Q)=z+c_\tau(Q)$ or $(**)$ for all $\tau\in D'$ we have
$T_\tau(Q)=-z+c_\tau(Q)$.
Without loss of generality, we may assume the case $(**)$,
that is, for all $\tau \in D'$ we have $T_\tau(Q)=-z+c_\tau(Q)$.
By Theorem \ref{Thm: algebraic lifting for torus},
$F\colon\mathcal{C}^n_\pi(V(1,1))\to
\mathcal{C}^n_\pi(V(1,1))$
can be lifted to a strictly equivariant fibred
morphism $f\colon\mathcal{E}^n_\pi(V(1,1))\to
\mathcal{E}^n_\pi(V(1,1))$ commuting with $\mathbf{S}(n)$-action 
(see the end of the proof of Theorem \ref{Thm: tame torus}).
The condition $(**)$ implies that the fibred morphism 
$h=\Id+f\colon\mathcal{E}^n_\pi(V(1,1))
\to(V(1,1))^n_\pi$ maps the set
${\theta}^{-1}(D')\subset \mathcal{E}^n_\pi(V(1,1))$
into the closed complex subspace $S\subset(V(1,1))^n_\pi$ (see
Definition \ref{Def: fibred configuration space}).
The morphism $h$ is continuous and ${\theta}^{-1}(D')$ is dense in
${\theta}^{-1}(D)$, thus $h({\theta}^{-1}(D))\subseteq S$.
Since $\mathcal{E}^n_\pi(V(1,1))$ is irreducible 
(see Remark \ref{Rmk: fibred configuration spaces are irreducible}),
it follows that $h(\mathcal{E}^n_\pi(V(1,1)))\subseteq S
\cong V(1,1)$. By the definition of $h$, it commutes with 
$\mathbf{S}(n)$-action, thus, $h$ is $\mathbf{S}(n)$-invariant
and it induces a fibred morphism 
$g\colon\mathcal{C}^n_\pi(V(1,1))\to V(1,1)$
with the desired properties.
\end{proof}

\begin{Remark}
If $F$ is an automorphism, the statement of the above
theorem holds true for $n =3, 4$.
\end{Remark}


\section{Acknowledgements}
\noindent This paper is based upon my PhD
thesis supported by Technion, Israel Institute of Technology.
I wish to thank V. Lin who
introduced me to the problem and encouraged me to work on it.


\bibliographystyle{plane}


\end{document}

%% file: cd.tex
\catcode`\@=11

 \def\AMSTeXfeatures{\Plainheads 
   \let\current@vert=\AMS@vert}

 \def\Plainheads{\sh@ftdiam=0.05em
   \getlabeldims
   \let\vshaftfill=\plnvsolidfill
   \let\hshaftfill=\plnhsolidfill
   \let\th@rhead=\plnrhead
   \let\th@lhead=\plnlhead
   \let\th@dnhead=\plndnhead
   \let\th@uphead=\plnuphead}
 
 \def\glet{\global\let}

 \def\LaTeXfeatures{\catcode`\@=11
   \ifx\@clnwd\undefined \nol@g
      \input ltxcode.tex \dol@g \fi
   \ltxheads \let\current@vert=\new@vert
   \providelto \catcode`\@=\active}

 \def\nol@g{\def\wlog{\edef\garbage}}
 \def\dol@g{\let\wlog=\wl@g} \let\wl@g=\wlog
 \nol@g 

 \newbox\ltobox
 \def\providelto{{\setbox\z@=
   \hbox{$\to$}\minharrlen=\wd\z@
   \global\setbox\ltobox=\hbox{$\activeat>>>$}}
   \def\lto{\mathrel{\copy\ltobox}}}

 \def\ltxheads{\sh@ftdiam=\@wholewidth
   \getlabeldims
   \let\vshaftfill= \ltxvsolidfill
   \let\hshaftfill=\ltxhsolidfill
   \let\th@rhead=\ltxrhead
   \let\th@lhead=\ltxlhead
   \let\th@dnhead=\ltxdnhead
   \let\th@uphead=\ltxuphead}
 {\catcode`\@=\active
   \gdef@#1{\csname #1\string@at\endcsname}
   \glet\activeat=@}
 \def\def@#1{\expandafter\def\csname #1@at\endcsname}

 \def@>#1>#2>{\@rrow R{#1}{#2}}
 \def@<#1<#2<{\@rrow L{#1}{#2}}
 \def@ V#1V#2V{\@rrow V{#1}{#2}}
 \def@ A#1A#2A{\@rrow A{#1}{#2}}
 \def@/#1/#2/#3/{\@rrow{#1}{#2}{#3}}
 \def@.{\ifodd\row\ifmmode\noharrow
     \else\leavevmode.\spacefactor3000 \fi
   \else\novarrow\fi}
 \def@={\ifodd\row\harrow\hequalfill{}{}%
   \else\varrow\vequalfill{}{}\fi}
 \def@:#1{\ifx=#1\harrow\deffill{}{}%
   \else\leavevmode\null:#1\fi}
 \def@|{\current@vert}
  \def\AMS@vert{\varrow\vequalfill{}{}}
  \def\new@vert#1|#2|{\ifodd\row
   \let\nextarrow\vertexvarrow
   \else\let\nextarrow\varrow\fi
   \nextarrow\vshaftfill{#1}{#2}}
 \def@-{\ifmmode\let\next\hl@ne
   \else\let\next\AMSatdash \fi \next}
  \def\hl@ne#1-#2-{\harrow\hshaftfill{#1}{#2}}
  \def\AMSatdash{\let\next\relax\leavevmode
    \def\next@{\ifx\next-%
      \def\next-{\futurelet\next\nextii@}%
     \else\def\next{\hbox{-}}\fi\next}%
    \def\nextii@{\ifx\next-\def\next-{\hbox{---}}%
      \else\def\next{\hbox{--}}\fi\next}%
    \futurelet\next\next@}
 \def@(#1){\tweenarrows{#1}}
 \def@[#1]{\setsp@n#1\relax\activeat}
 \def\fiberbox{\hbox{$\vcenter{\hr@le\hbox{\vr@le
   \kern1ex\vbox{\kern1.2ex}\vr@le}\hr@le}$}}
  \def\hr@le{\hrule height \sh@ftdiam}
  \def\vr@le{\vrule width \sh@ftdiam}
 \def@+#1+#2+#3+{\ifodd\row \harrow{#1}{#2}{#3}%
   \else \varrow{#1}{#2}{#3}\fi}

 \def\Dnarrfill{\vequalfill\Dnhe@d}
 \def\Uparrfill{\Uphe@d\vequalfill}
 
 \def\ontofill{\rtarrfill\kern-0.3em 
   \th@rhead\kern 0.3em} 

 \def\rtarrfill{\hshaftfill\th@rhead}
 \def\ltarrfill{\th@lhead\hshaftfill}
 \def\dnarrfill{\vshaftfill\th@dnhead}
 \def\uparrfill{\th@uphead\vshaftfill}
 \def\hequalfill{\plnhfill=}
 \def\deffill{:\plnhfill=}
 \def\plnvextfill#1{\setbox\z@
   \hbox{\the\textfont3 #1}%
   \dimen@=\dp\z@\advance\dimen@\ht\z@
   \copy\z@ \kern-\dimen@ 
   \cleaders\copy\z@ \vfill
   \kern-\dimen@ 
   \box\z@}
 \def\plnhfill#1{$\m@th\mkern-1.5mu\mathord#1\mkern-6mu
    \cleaders\hbox{$\mkern-2mu\mathord#1\mkern-2mu$}\hfill
    \mkern-6mu\mathord#1\mkern-1.5mu$}
 \def\vequalfill{\plnvextfill{\char'167}}
 \def\plnvsolidfill{\plnvextfill{\char'077}}
 \def\plnhsolidfill{\plnhfill-}
 \def\ltxhsolidfill{\leaders\hrule height\topofshaft depth\botofshaft
   \hfill}
 \def\ltxvsolidfill{\leaders\vrule width\sh@ftdiam\vfill}
 \def\hdashfill{\hd@sh\wd@sh
   \xleaders \hbox{\wd@sh\hd@sh\wd@sh}\hfill
   \wd@sh\hd@sh}
 \def\vdashfill{\vd@sh\wd@sh
   \xleaders \vbox{\wd@sh\vd@sh\wd@sh}\vfill
   \wd@sh\vd@sh}
 \def\dashed{\ifinmeasureCD\else
    \ifodd\row\option{\let\hshaftfill=\hdashfill}%
   \else\option{\let\vshaftfill=\vdashfill}\fi\fi}

 \newdimen\CDstrutht  \newdimen\CDstrutdp
 \newdimen\CDstrutlen \CDstrutlen=\CDstrutht
   \advance\CDstrutlen by \CDstrutdp

 \def\CDstrut{\vrule
   height \ifnum\row=1 \z@\else\CDstrutht \fi
   depth \ifnum\row=\numrows \z@ \else\CDstrutdp \fi
   width\z@}

 \newdimen\CDarrsurr \CDarrsurr=0.375em
 \newdimen\CDdashlen
 \newdimen\CDvarrlen \CDvarrlen=1.5\baselineskip
 \newdimen\minharrlen 
  \setbox\z@\hbox{$\longrightarrow$} \minharrlen=\wd\z@
 \newdimen\minCDharrlen \minCDharrlen=2.5em 
\newdimen \minc@lwd
\def\findminc@lwd{\minc@lwd=2\CDarrsurr
  \advance\minc@lwd\minCDharrlen}

 \newdimen\sh@ftdiam

 \newdimen\labelsurr \labelsurr=1.25 em

\newcount\sp@ncnt \sp@ncnt=\@ne
\newcount\sp@ncnt@ \sp@ncnt@=\@ne
\newdimen\@rrwd \newdimen\@rrdp

 \def\adjustbot#1{\option{\advance\@rrdp#1\relax}}
 
\def\pushvertex#1{\global\p@shlen#1\relax
   \global\let\maybepush=\dopush}

 \newdimen\p@shlen \p@shlen=\z@

 \let\maybepush=\relax
 \def\dopush{\ifinmeasureCD 
   \advance\locdimen by -\p@shlen 
   \else\advance \@rrwd by -\p@shlen \fi 
   \global\let\maybepush=\relax \global\p@shlen=\z@\relax}

 \def\span@ne{\global\sp@ncnt=\@ne\relax}
 \def\setsp@n#1#2{\global\sp@ncnt=#1\relax
   \ifx\relax#2\relax\else\global\sp@ncnt@=#2\relax\fi}

 \def\plnrhead{\llap{$\rightarrow\mkern-1.5mu$}}
 \def\plnlhead{\rlap{$\mkern-1.5mu\leftarrow$}}

 \def\clap#1{\hbox to \z@{\hss #1\hss}}

 \def\plndnhead{\hbox{\the\textfont3 \char'171}}
 \def\plnuphead{\hbox{\the\textfont3 \char'170}}
 \def\Dnhe@d{\hbox{\the\textfont3 \char'177}}
 \def\Uphe@d{\hbox{\the\textfont3 \char'176}}

 \def\ltxrhead{\raise\@xisheight
   \llap{\smash{\@linefnt\@getrarrow(1,0)}}}
 \def\ltxlhead{\raise\@xisheight
   \rlap{\@linefnt\@getlarrow(-1,0)}}
 \def\ltxuphead{\setbox\z@=\rlap{%
   \kern\@halfwidth\@linefnt\char'66}%
   \copy\z@\kern-\ht\z@}
 \def\ltxdnhead{\setbox\z@=\rlap{%
   \kern\@halfwidth\@linefnt\char'77}%
   \ht\z@=\z@\box\z@}

 \def\wd@sh{\kern0.5\CDdashlen}
 \def\hd@sh{\vrule height\topofshaft depth\botofshaft
    width\CDdashlen}
 \def\vd@sh{\hrule height\CDdashlen
   depth\z@ width\sh@ftdiam}

\def\xylist{14{3434}13{2414}12{1723}%
  23{1413}34{1153}11{0867}43{0707}%
  32{0580}21{0414}31{0291}41{0}}
\newcount\tgtcnt@
\def\find@xyargs{\dimen@=\@rrdp
  \advance\dimen@ by \CDstrutlen
  \tgtcnt@=\dimen@ \dimen@=\@rrwd 
  \divide\dimen@ by \@m 
  \divide \tgtcnt@ by \dimen@ 
  \expandafter\testxy\xylist\relax
  \unitlength=\@xarg\@rrdp
  \divide\unitlength by\@yarg\relax}
\def\testxy#1#2#3{\ifnum\tgtcnt@>#3
    \@xarg=#1\relax \@yarg=#2\relax
    \let\next=\ignorerest
  \else\let\next\testxy\fi\next}
\def\ignorerest#1\relax{\relax}

\let\scalefactor=\@ne
\def\SWarrow{\find@xyargs\vector
  (-\@xarg,-\@yarg)\scalefactor\hskip-\wd\@linechar}
\def\NWarrow{\find@xyargs\vector
  (-\@xarg,\@yarg)\scalefactor\hskip-\wd\@linechar}
\def\NEarrow{\find@xyargs\vector
  (\@xarg,\@yarg)\scalefactor}
\def\SEarrow{\find@xyargs\vector
  (\@xarg,-\@yarg)\scalefactor}
\def\rightupline{\find@xyargs\@linelen=\scalefactor
     \unitlength\@sline}
\def\rightdownline{\find@xyargs\@yarg=-\@yarg\relax
     \@linelen=\scalefactor\unitlength\@sline}

\def\Sim{\ifodd\row\setbox\z@=\hbox{$\sim$}\dimen@=\ht\z@
 \advance\dimen@ by -\@xisheight
  \vbox{\box\z@\kern-\@xisheight\kern\dimen@}%
  \else\hbox{$\wr$}\fi}

\def\harrow#1#2#3{\inmeasureCDtrue\findminarrwd
  {#2}{#3}{\sp@ncnt\minharrlen}\inmeasureCDfalse\span@ne
  \mathrel{\hbox{\options\hplace{#1}\ulabel{#2}\dlabel{#3}}}}

\def\noharrow{\harrow\hfill{}{}}
\def\vertexvarrow#1#2#3{\findarrdp \@rrwd=\z@ \setsp@n\@ne\@ne
  \vbox to \z@{\kern-1.2\CDstrutht
  \rlap{\options\vplace{#1}\llabel{#2}\rlabel{#3}}\vss}}

\newif\ifinmeasureCD
\def\measurelabel#1{\setbox\z@
  \hbox{$\scriptstyle#1\kern\labelsurr$}%
  \ifdim\wd\z@>\@rrwd \@rrwd=\wd\z@\fi}
\def\findminarrwd#1#2#3{\@rrwd=#3\relax
   \measurelabel{#1}\measurelabel{#2}}
\def\findCDarrwd#1#2{\@rrwd=\minCDharrlen
   \measurelabel{#1}\measurelabel{#2}%
  }

\newcount\row \row=\@ne \newcount\col \col=\@ne 
 \newcount\numrows 
\numrows=\@ne
 \newcount\numcols
\newcount\arrspan \newdimen\vrtxhalfwd  \newbox\tempbox

\def\DANABUG{\advance\col by \@ne
 \@rrwd=\minCDharrlen
  \advance\@rrwd by \vrtxhalfwd
  \advance\@rrwd by \CDarrsurr
  \ifnum\col>\numcols \numcols=\col
     \newlocdimen{col\the\col}\locdimen=\@rrwd 
  \else \ifdim\@rrwd>\c@l \c@l=\@rrwd\fi\fi}

\def\drop#1\\{
  \findvrtxhalfsum\DANABUG\advance\row by 2 \measureinit}

\def\measureinit{\col=\@ne \vrtxhalfwd=-\CDarrsurr\arrspan=\@ne\@rrwd=\z@
   \setbox\tempbox=\hbox\bgroup$}
\def\measure{
  \let\harrow\measureCDarrow
  \let\CDCR=\measureCR 
   \findminc@lwd 
  \inmeasureCDtrue
  \row=\@ne \numcols=\z@ \measureinit}

\def\endmeasure{\findvrtxhalfsum\DANABUG
  \numrows=\row 
  \inmeasureCDfalse}

\def\newlocdimen#1{\advance\dimenc@unt by \@ne
  \ifnum\dimenc@unt<\insc@unt
     \else\errmessage{No room for the CD}\fi
  \dimendef\locdimen=\dimenc@unt
  \expandafter\dimendef\csname#1\endcsname=\dimenc@unt}

 \def\r@wc@l{\csname row\the\row col\the\col\endcsname}
 \def\c@l{\csname col\the\col\endcsname}

 \def\findvrtxhalfsum{$\egroup
  \newlocdimen{row\the\row col\the\col}
  \locdimen=\vrtxhalfwd 
  \vrtxhalfwd=0.5\wd\tempbox 
  \advance\vrtxhalfwd by \CDarrsurr
  \advance\locdimen by \vrtxhalfwd 
  \advance\@rrwd by \locdimen 
  \maybepush
  \divide\@rrwd by \arrspan\relax
  \ifdim\@rrwd<\minc@lwd
    \ifnum\col>\@ne \@rrwd=\minc@lwd\fi \fi
  \loop 
    \ifnum\col>\numcols \numcols=\col
       \newlocdimen{col\the\col}
       \locdimen=\@rrwd 
    \else \ifdim\@rrwd>\c@l \c@l=\@rrwd\fi \fi
   \ifnum\arrspan>\@ne
      \advance\arrspan by -1 \advance\col by \@ne
  \repeat }

 \def\measureCDarrow#1#2#3{\findvrtxhalfsum
   \arrspan=\sp@ncnt\relax\global\sp@ncnt=1\relax
   \advance\col by \@ne
   \findCDarrwd{#2}{#3}%
   \setbox\tempbox=\hbox\bgroup$}

 \newcount\dr@tn \dr@tn=\z@
 \def\locate#1:#2{\ifinmeasureCD\else
   \count@=-#1
   \multiply\count@ by 2
   \advance\count@ by #2
   \dimen@=\count@\@rrwd
   \ifnum\dr@tn=\@ne\relax \else\dimen@=-\dimen@ \fi
   \dimen@i=\@rrdp
   \ifnum\dr@tn>\z@\advance\dimen@i by \CDstrutlen \fi
   \dimen@i=\count@\dimen@i
   \count@=#2 \multiply\count@ by 2
   \divide\dimen@ by \count@
   \divide\dimen@i by \count@
   \lift\dimen@i\nudge\dimen@\fi}

\def\betweenCDrows{\advance\row by \@ne \col=\@ne
\options}

\def\hbegin{\hbox\bgroup\kern\c@l \kern-\r@wc@l$}
\def\hend{$\glet\maybepush\relax \CDstrut\egroup}
\def\vbegin{\setbox\tempbox=\hbox\bgroup$}
\def\vend{$\egroup\ht\tempbox=\z@\dp\tempbox\CDvarrlen
  \box\tempbox}
\def\setCD{\let\harrow=\setCDarrow
  \let\CDCR=\setCR 
  \row=\@ne \col=\@ne \hbegin}
\let\endsetCD=\hend 

\def\findarrwd{\@rrwd=\z@ \count@=\col \advance\count@ by\sp@ncnt
  \loop\ifnum\count@>\col \advance\count@ by -1
      \advance\@rrwd by\csname col\the\count@\endcsname\repeat}
\def\setCDarrow#1#2#3{\kern\CDarrsurr\advance\col by \@ne
  \findarrwd \advance\@rrwd by -\r@wc@l  
  \@rrdp=\z@ 
  \maybepush
  \advance\col by -\@ne \advance\col by \sp@ncnt \span@ne
  \hbox to \@rrwd{\options
   \@rrwd=\scalefactor\@rrwd\hss
   \hplace{#1}\ulabel{#2}\dlabel{#3}\hss}%
   \kern\CDarrsurr}

\newdimen\labspacei 
\newdimen\labspaceii 

\newdimen\@xisheight
  \@xisheight=\the\fontdimen22\textfont2
\newdimen\labelskip
  \labelskip=\the\fontdimen10\textfont3 
\newdimen\topofshaft
\newdimen\botofshaft
\newdimen\botofulabel
\newdimen\topofdlabel
\def\getlabeldims{
  \topofshaft=0.5\sh@ftdiam
  \botofshaft=\topofshaft
  \advance\topofshaft by \@xisheight  
  \advance\botofshaft by -\@xisheight  
  \botofulabel=\topofshaft
  \advance\botofulabel by \labelskip
  \topofdlabel=\botofshaft
  \advance\topofdlabel by \labelskip}

\def\ulabel{\ifnum\row=\@ne\let\next\ulabeli
   \else\let\next\ulabellap\fi\next}
\def\ulabeli#1{\vbox{
  \clap{\kern-\@rrwd$\scriptstyle#1$}%
  \kern\botofulabel}\maybeoffset}
\def\ulabellap#1{\vbox to \z@{\vss
  \clap{\kern-\@rrwd$\scriptstyle#1$}%
  \kern\botofulabel}\maybeoffset}
\def\dlabel{\ifnum\row=\numrows\let\next\dlabeli
   \else\let\next\dlabellap\fi\next}
\def\dlabeli#1{\vtop{\kern\topofdlabel
  \clap{\kern-\@rrwd$\scriptstyle#1$}%
  }\maybeoffset}
\def\dlabellap#1{\vbox to \z@{\kern\topofdlabel
  \clap{\kern-\@rrwd$\scriptstyle#1$}%
  \vss}\maybeoffset}
\def\rlabel#1{\vbox to \z@{\vss
  \rlap{\kern\labelskip$\scriptstyle#1$}%
  \vss\kern-\@rrdp}\maybeoffset}
\def\llabel#1{\vbox to \z@{\vss
  \llap{$\scriptstyle#1$\kern\labelskip}%
  \vss\kern-\@rrdp}\maybeoffset}
\def\swlabel#1{\vtop{\kern0.5\@rrdp
  \llap{$\scriptstyle#1$\kern\labelskip\kern-0.5\@rrwd}
  }\maybeoffset}
\def\nwlabel#1{\vbox{
  \llap{$\scriptstyle#1$\kern\labelskip\kern-0.5\@rrwd}%
  \kern-0.5\@rrdp}\maybeoffset}
\def\selabel#1{\vtop{\kern0.5\@rrdp
  \rlap{\kern0.5\@rrwd\kern\labelskip$\scriptstyle#1$}%
  }\maybeoffset}
\def\nelabel#1{\vbox{
  \rlap{\kern0.5\@rrwd\kern\labelskip$\scriptstyle#1$}%
  \kern-0.5\@rrdp}\maybeoffset}
\def\cplace#1{\vbox to \z@{\vss
  \clap{$#1$\kern-\@rrwd}%
  \kern-\@rrdp\vss}\maybeoffset}
\def\hplace#1{\hbox to \@rrwd{#1}\maybeoffset}
\def\vplace#1{\clap{\vbox to \z@{#1\kern-\@rrdp}}\maybeoffset}

\newdimen\nudgeamount \nudgeamount=\z@
\newdimen\liftamount \liftamount=\z@
\let\maybeoffset\relax
\newbox\offsetbox \newdimen\lastheight
\def\dooffset{
  \setbox\offsetbox=\lastbox \lastheight=\ht\offsetbox 
  \setbox\offsetbox=\vbox{\kern-\liftamount\box\offsetbox}%
  \ht\offsetbox=\lastheight
  \kern\nudgeamount\box\offsetbox\kern-\nudgeamount
  \global\nudgeamount=\z@ \global\liftamount=\z@
  \glet\maybeoffset=\relax}
\def\nudge#1{\ifinmeasureCD\else
  \global\advance\nudgeamount#1\relax
  \global\let\maybeoffset\dooffset\fi}
\def\lift#1{\ifinmeasureCD\else
  \global\advance\liftamount#1\relax
  \global\let\maybeoffset\dooffset\fi}

\def\findarrdp{\@rrdp=\CDvarrlen
  \ifnum\sp@ncnt@>1
    \advance\@rrdp by \CDstrutlen
    \multiply\@rrdp by \sp@ncnt@
    \advance\@rrdp by -\CDstrutlen \fi
 }

\def\varrow#1#2#3{\ifnum\sp@ncnt>\@ne 
     \sp@ncnt@=\sp@ncnt\relax\fi
  \findarrdp \@rrwd=\z@ 
  \kern\c@l
   \hbox to \z@{\options
   \@rrdp=\scalefactor\@rrdp
    \hss\vplace{#1}\llabel{#2}\rlabel{#3}\hss}%
  \global\advance\col by \@ne \setsp@n\@ne\@ne
  }

\def\novarrow{\varrow\vfill{}{}}

\def\tweenarrows#1{\findarrwd \findarrdp \setsp@n\@ne\@ne
  \rlap{\options\cplace{#1}}}

\def\usarrow #1#2#3{\dr@tn=\@ne
  \findarrwd \findarrdp \setsp@n\@ne\@ne 
  \rlap{\options\cplace{#1}\nwlabel{#2}\selabel{#3}}%
  \dr@tn=\z@}
\def\dsarrow #1#2#3{\dr@tn=\tw@
  \findarrwd \findarrdp \setsp@n\@ne\@ne 
  \rlap{\options\cplace{#1}\swlabel{#2}\nelabel{#3}}%
  \dr@tn=\z@}
 \def\@rrow#1{\csname #1@rrow\endcsname}
 \def\R@rrow{\harrow \rtarrfill}
 \def\L@rrow{\harrow \ltarrfill}
 \def\V@rrow{\varrow \dnarrfill}
 \def\A@rrow{\varrow \uparrfill}
 \def\SE@rrow{\dsarrow \SEarrow}
 \def\NW@rrow{\dsarrow \NWarrow}
 \def\SW@rrow{\usarrow \SWarrow}
 \def\NE@rrow{\usarrow \NEarrow}
 \def\DS@rrow{\dsarrow \dnslope}
 \def\US@rrow{\usarrow \upslope}
 \def\upslope{\find@xyargs
       \@linelen=\unitlength\@sline}
 \def\dnslope{\find@xyargs\@yarg=-\@yarg\relax
       \@linelen=\unitlength\@sline}

\newtoks\optionlist 
\optionlist={}
\let\options\relax
\def\dooptions{\the\optionlist\global\optionlist={}%
  \glet\options=\relax}
\def\option#1{\ifinmeasureCD\else
  \glet\options=\dooptions
  \global\optionlist=\expandafter{\the\optionlist\relax#1}\fi}
\def\wider#1{\ifinmeasureCD\else
   \option{\advance\@rrwd by #1}\fi}
\def\deeper#1{\ifinmeasureCD\else
   \option{\advance\@rrdp by #1}\fi}
\def\arrowscale#1{\ifinmeasureCD\else
   \option{\def\scalefactor{#1}}\fi}

{\def\\{\global\let\sptoken= }\\ }

\def\CR{\futurelet\nexttok\testCR}
\def\testCR{\ifx\nexttok\sptoken
   \let\next\eatspaceCR\else\let\next\CDCR\fi\next}
\def\eatspaceCR#1 {\CR}
\def\measureCR{\ifx\nexttok\endmeasure\let\nextCR\relax
    \else\let\nextCR\drop\fi\nextCR}
\def\setCR{\ifodd\row
  \ifx\nexttok\endsetCD\else\hend\betweenCDrows\vbegin\fi
  \else\vend\betweenCDrows\hbegin\fi}

\countdef\dimenc@unt=11
\def\CD#1\endCD{
   \begingroup\let\\=\CR
  \m@th\offinterlineskip
   \measure#1\endmeasure\null\,\vcenter{\setCD#1\endsetCD}\,
   \endgroup
    }

\ifx\@clnwd\undefined \nol@g\else\catcode`\ =14\relax\fi
 \font\@linefnt=line10 
 \newcount\@tempcnta
 \newcount\@tempcntb
 \newdimen\@tempdima
 \newdimen\@tempdimb
 \newdimen\@wholewidth
 \newdimen\@halfwidth
   \@wholewidth\fontdimen8\@linefnt \@halfwidth .5\@wholewidth
 \newdimen\unitlength
 \newcount\@xarg
 \newcount\@yarg
 \newcount\@yyarg
 \newbox\@linechar
 \newdimen\@linelen
 \newdimen\@clnwd
 \newdimen\@clnht
 \newif\if@negarg
 
 \def\@whilenoop#1{}

 \def\@whiledim#1\do #2{\ifdim #1\relax#2\@iwhiledim{#1\relax#2}\fi}

 \def\@iwhiledim#1{\ifdim #1\let\@nextwhile=\@iwhiledim 
         \else\let\@nextwhile=\@whilenoop\fi\@nextwhile{#1}}

 \def\@sline{\ifnum\@xarg< 0 \@negargtrue \@xarg -\@xarg \@yyarg -\@yarg
   \else \@negargfalse \@yyarg \@yarg \fi
 \ifnum \@yyarg >0 \@tempcnta\@yyarg \else \@tempcnta -\@yyarg \fi
 \ifnum\@tempcnta>6 \@badlinearg\@tempcnta0 \fi
 \ifnum\@xarg>6 \@badlinearg\@xarg 1 \fi
 \setbox\@linechar\hbox{\@linefnt\@getlinechar(\@xarg,\@yyarg)}%
 \ifnum \@yarg >0 \let\@upordown\raise \@clnht\z@
    \else\let\@upordown\lower \@clnht \ht\@linechar\fi
 \@clnwd=\wd\@linechar
 \if@negarg \hskip -\wd\@linechar \def\@tempa{\hskip -2\wd\@linechar}\else
      \let\@tempa\relax \fi
 \@whiledim \@clnwd <\@linelen \do
   {\@upordown\@clnht\copy\@linechar
    \@tempa
    \advance\@clnht \ht\@linechar
    \advance\@clnwd \wd\@linechar}%
 \advance\@clnht -\ht\@linechar
 \advance\@clnwd -\wd\@linechar
 \@tempdima\@linelen\advance\@tempdima -\@clnwd
 \@tempdimb\@tempdima\advance\@tempdimb -\wd\@linechar
 \if@negarg \hskip -\@tempdimb \else \hskip \@tempdimb \fi
 \multiply\@tempdima \@m
 \@tempcnta \@tempdima \@tempdima \wd\@linechar \divide\@tempcnta \@tempdima
 \@tempdima \ht\@linechar \multiply\@tempdima \@tempcnta
 \divide\@tempdima \@m
 \advance\@clnht \@tempdima
 \ifdim \@linelen <\wd\@linechar
    \hskip \wd\@linechar
   \else\@upordown\@clnht\copy\@linechar\fi}
 
 \def\@getlinechar(#1,#2){\@tempcnta#1\relax\multiply\@tempcnta 8
 \advance\@tempcnta -9 \ifnum #2>0 \advance\@tempcnta #2\relax\else
 \advance\@tempcnta -#2\relax\advance\@tempcnta 64 \fi
 \char\@tempcnta}
 
 \def\vector(#1,#2)#3{\@xarg #1\relax \@yarg #2\relax
 \@tempcnta \ifnum\@xarg<0 -\@xarg\else\@xarg\fi
 \ifnum\@tempcnta<5\relax
 \@linelen=#3\unitlength
 \ifnum\@xarg =0 \@vvector 
   \else \ifnum\@yarg =0 \@hvector \else \@svector\fi
 \fi
 \else\@badlinearg\fi}
 
 \def\@svector{\@sline
 \@tempcnta\@yarg \ifnum\@tempcnta <0 \@tempcnta=-\@tempcnta\fi
 \ifnum\@tempcnta <5
   \hskip -\wd\@linechar
   \@upordown\@clnht \hbox{\@linefnt  \if@negarg 
   \@getlarrow(\@xarg,\@yyarg) \else \@getrarrow(\@xarg,\@yyarg) \fi}%
 \else\@badlinearg\fi}
 
 \def\@getlarrow(#1,#2){\ifnum #2 =\z@ \@tempcnta='33\else
 \@tempcnta=#1\relax\multiply\@tempcnta \sixt@@n \advance\@tempcnta
 -9 \@tempcntb=#2\relax\multiply\@tempcntb \tw@
 \ifnum \@tempcntb >0 \advance\@tempcnta \@tempcntb\relax
 \else\advance\@tempcnta -\@tempcntb\advance\@tempcnta 64
 \fi\fi\char\@tempcnta}
 
 \def\@getrarrow(#1,#2){\@tempcntb=#2\relax
 \ifnum\@tempcntb < 0 \@tempcntb=-\@tempcntb\relax\fi
 \ifcase \@tempcntb\relax \@tempcnta='55 \or 
 \ifnum #1<3 \@tempcnta=#1\relax\multiply\@tempcnta
 24 \advance\@tempcnta -6 \else \ifnum #1=3 \@tempcnta=49
 \else\@tempcnta=58 \fi\fi\or 
 \ifnum #1<3 \@tempcnta=#1\relax\multiply\@tempcnta
 24 \advance\@tempcnta -3 \else \@tempcnta=51\fi\or 
 \@tempcnta=#1\relax\multiply\@tempcnta
 \sixt@@n \advance\@tempcnta -\tw@ \else
 \@tempcnta=#1\relax\multiply\@tempcnta
 \sixt@@n \advance\@tempcnta 7 \fi\ifnum #2<0 \advance\@tempcnta 64 \fi
 \char\@tempcnta}
\catcode`\ =10

\dol@g 
\catcode`\@=\active
\LaTeXfeatures